\documentclass[leqno]{article}

\usepackage[frenchb,english]{babel}
\usepackage[T1]{fontenc}

\usepackage{amsmath}
\usepackage{amssymb}
\usepackage{amsfonts}
\usepackage{enumerate}
\usepackage{vmargin}
\usepackage[all]{xy}
\usepackage{mathrsfs}

\setmarginsrb{3cm}{3cm}{2.5cm}{2cm}{0cm}{0cm}{1.5cm}{4cm}
\newcommand{\Mat}{\ensuremath{\mathbb{M}}}
\newcommand{\Z}{\ensuremath{\mathbb{Z}}}
\newcommand{\R}{\ensuremath{\mathbb{R}}}
\newcommand{\C}{\ensuremath{\mathbb{C}}}
\newcommand{\T}{\ensuremath{\mathbb{T}}}

\newcommand{\tr}{\ensuremath{\mathop{\rm Tr\,}\nolimits}}

\newcommand{\Id}{{\rm Id}}
\newcommand{\VN}{{\rm VN}}

\renewcommand{\leq}{\ensuremath{\leqslant}}
\renewcommand{\geq}{\ensuremath{\geqslant}}
\newcommand{\n}{\noindent}

\newcommand{\qed}{\hfill \vrule height6pt  width6pt depth0pt}

\newcommand{\ull}{\mathcal{U}}

\newcommand\be{\begin{eqnarray}}
\newcommand\ee{\end{eqnarray}}
\newcommand{\norm}[1]{ \| #1  \|}
\newcommand{\bnorm}[1]{ \big\| #1  \big\|}

\newcommand{\Bgnorm}[1]{ \Bigg\| #1  \Bigg\|}
\newcommand{\xra}{\xrightarrow}

\newcommand{\ot}{\otimes}

\def\Aut{{\rm Aut} }

\def\card{{\rm card} }
\newcommand{\ovl}{\overline}
\newcommand{\otvn}{\ovl\ot}

\newtheorem{thm}{Theorem}[section]
\newtheorem{defi}[thm]{Definition}
\newtheorem{prop}[thm]{Proposition}
\newtheorem{conj}[thm]{Conjecture}
\newtheorem{quest}[thm]{Question}
\newtheorem{Property}[thm]{Property}
\newtheorem{cor}[thm]{Corollary}
\newtheorem{lemma}[thm]{Lemma}

\newtheorem{remark}[thm]{Remark}

\newenvironment{preuve}[1][]{\noindent {\it Proof #1} : }{\hbox{~}\qed
\smallskip
}

\newcommand{\beq}{\begin{equation}}
\newcommand{\eeq}{\end{equation}}

\numberwithin{equation}{section}

\begin{document}
\selectlanguage{english}
\title{\bfseries{On Matsaev's conjecture for contractions on noncommutative $L^p$-spaces}}
\date{}
\author{\bfseries{C\'edric Arhancet}}

\maketitle

%%%%%%%%%%%%%%%%%%%%%%%%%%%%%%%%%%%%%%%%%%%%%%%%%%%%%%%%%%%%%%%%
%%%%%%%%%%%%%%%%%%%%%%%%%%%%%%%%%%%%%%%%%%%%%%%%%%%%%%%%%%%%%%%%
\begin{abstract}
We exhibit large classes of contractions on noncommutative
$L^p$-spaces which satisfy the noncommutative analogue of Matsaev's
conjecture, introduced by Peller, in 1985. In particular, we prove
that every Schur multiplier on a Schatten space $S^p$ induced by a
contractive Schur multiplier on $B(\ell^2)$ associated with a real
matrix satisfy this conjecture. Moreover, we deal with analogue
questions for $C_0$-semigroups. Finally, we disprove a conjecture of
Peller concerning norms on the space of complex polynomials arising
from Matsaev's conjecture and Peller's problem. Indeed, if $S$
denotes the shift on $\ell^p$ and $\sigma$ the shift on the Schatten
space $S^p$, the norms $\bnorm{P(S)}_{\ell^p \xra{}\ell^p}$ and
$\bnorm{P(\sigma)\ot \Id_{S^p}}_{S^p(S^p) \xra{}S^p(S^p)}$  can be
different for a complex polynomial $P$.
\end{abstract}

% and every Fourier multiplier on
%$L_p\big(\VN(G)\big)$ induced by a contractive Fourier multiplier on
%$VN(G)$ associated with a real function, where $\VN(G)$ is the von
%Neumann algebra of an amenable discrete group $G$,

%%%%%%%%%%%%%%%%%%%%%%%%%%%%%%%%%%%%%%%%%%%%%%%%%%%%%%%%%%%%%%%%
%%%%%%%%%%%%%%%%%%%%%%%%%%%%%%%%%%%%%%%%%%%%%%%%%%%%%%%%%%%%%%%%

\makeatletter
 \renewcommand{\@makefntext}[1]{#1}
 \makeatother \footnotetext{\noindent
 This work is partially supported by ANR 06-BLAN-0015.
 }

%%%%%%%%%%%%%%%%%%%%%%%%%%%%%%%%%%%%%%%%%%%%%%%%%%%%%%%%%%%%%%%%
%%%%%%%%%%%%%%%%%%%%%%%%%%%%%%%%%%%%%%%%%%%%%%%%%%%%%%%%%%%%%%%%

\section{Introduction}

%%%%%%%%%%%%%%%%%%%%%%%%%%%%%%%%%%%%%%%%%%%%%%%%%%%%%%%%%%%%%%%%
%%%%%%%%%%%%%%%%%%%%%%%%%%%%%%%%%%%%%%%%%%%%%%%%%%%%%%%%%%%%%%%%
To estimate the norms of functions of operators is an essential task
in Operator Theory. In this subject, V. V. Matsaev stated the
following conjecture in 1971, see \cite{Nik1}. For any $1\leq
p\leq\infty$, let $\ell^p\xra{S}\ell^p$ denote the right shift
operator defined by $S(a_0,a_1,a_2,\ldots)=(0,a_0,a_1,a_2,\ldots)$.

\begin{conj}
Suppose  $1<p< \infty$, $p\not=2$. Let $\Omega$ be a measure space
and let $L^p(\Omega)\xra{T}L^p(\Omega)$ be a contraction. For any
complex polynomial $P$, we have
\begin{eqnarray}
\label{Matsaev} \bnorm{P(T)}_{L^p(\Omega) \xra{}L^p(\Omega)}\leq
\bnorm{P(S)}_{\ell^p \xra{}\ell^p}.
\end{eqnarray}
\end{conj}

It is easy to see that (\ref{Matsaev}) holds true for $p=1$ and
$p=\infty$. Moreover, by using the Fourier transform, it is clear
that for $p=2$, (\ref{Matsaev}) is a consequence of von Neumann's
inequality. Finally, very recently and after the writing of this
paper, S. W. Drury \cite{Dr} found a counterexample in the case
$p=4$ by using computer.

For all other values of $p$, the validity of (\ref{Matsaev}) for any
contraction is open. It is well-known that (\ref{Matsaev}) holds
true for any positive contraction, more generally for all operators
$L^p(\Omega)\xra{T} L^p(\Omega)$ which admit a contractive majorant
$\big($i.e. there exists a positive contraction $\tilde{T}$
satisfying $|T(f)|\leq \tilde{T}(|f|)\big)$. This follows from the
fact that these operators admit an isometric dilation. We refer the
reader to  \cite{ALM}, \cite{AkS}, \cite{CoW}, \cite{Kit},
\cite{Nik2} and \cite{Pel1} for information and historical
background on this question.

In 1985, V.V. Peller \cite{Pel2} introduced a noncommutative version
of Matsaev's conjecture for Schatten spaces $S^p=S^p(\ell^2)$.
Recall that elements of $S^p$ can be regarded as infinite matrices
indexed by $\mathbb{N}\times \mathbb{N}$. Thus we define the linear
map $S^p \xra{\sigma}S^p$ as the shift `from NW to SE' which maps
any matrix
\begin{equation}\label{Shift}
\left[
  \begin{array}{cccc}
    a_{00} & a_{01} & a_{02} & \cdots\\
     a_{10}& a_{11} & a_{12} & \cdots \\
    a_{20} & a_{21}& a_{22} & \cdots\\
    \cdots & \cdots & \cdots & \cdots \\
  \end{array}
\right] \ \ \text{to} \ \ \left[
  \begin{array}{cccc}
    0 & 0 & 0 & \cdots \\
    0 & a_{00} & a_{01} & \cdots \\
    0 &  a_{10} &  a_{11} & \cdots \\
    \cdots & \cdots & \cdots & \cdots \\
  \end{array}
\right].
\end{equation}
Let $S^p(S^p)$ be the space of all matrices $[a_{ij}]_{i,j\geq 0}$
with entries $a_{ij}$ in $S^p$, which represent an element of the
bigger Schatten space $S^p(\ell^2\ot_2\ell^2)$. The algebraic tensor
product $S^p \ot S^p$ can be regarded as a dense subspace of
$S^p(S^p)$ in a natural way. Then the  mapping on $S^p(S^p)$ given
by (\ref{Shift}) is an isometry, which is the unique extension of
$\sigma \ot I_{S^p}$ to the space $S^p(S^p)$. (See Section 2 below
for more details on these matricial representations.) Peller's
question is as follows.

\begin{quest}
Suppose $1<p<\infty$, $p\not=2$. Let $S^p \xra{T}S^p$  be a
contraction on the Schatten space $S^p$. Do we have
\begin{equation}
\label{NCMatsaev} \bnorm{P(T)}_{S^p \xra{}S^p} \leq
\bnorm{P(\sigma)\ot Id_{S^p}}_{S^p(S^p) \xra{} S^p(S^p)}
\end{equation}
for any complex polynomial $P$?
\end{quest}

Peller observed that (\ref{NCMatsaev}) holds true when $T$ is an
isometry or when $S^p \xra{T}S^p $ is defined by $T(x)=axb$, where
$\ell^2\xra{a} \ell^2$ and $\ell^2\xra{b}\ell^2$ are contractions.

The Schatten spaces $S^p$ are basic examples of noncommutative
$L^p$-spaces. It is then natural to extend Peller's problem to this
wider context. This leads to the following question.

\begin{quest}
\label{quest matsaev semifinite} Suppose $1< p< \infty$, $p\not=2$.
Let $M$ be a semifinite von Neumann algebra and let $L^p(M)$ be the
associated noncommutative $L^p$-space. Let $L^p(M)\xra{T}L^p(M)$ be
a contraction. Do we have
\begin{eqnarray}
\label{LpMatsaev}
% \nonumber to remove numbering (before each equation)
\bnorm{P(T)}_{L^p(M) \xra{} L^p(M)} &\leq& \bnorm{P(\sigma)\ot
Id_{S^p}}_{S^p(S^p) \xra{}S^p(S^p)}
\end{eqnarray}
for any complex polynomial $P$?
\end{quest}

As in the commutative case, it is easy to see that (\ref{LpMatsaev})
holds true when $p=1$, $p=2$ or $p=\infty$. The main purpose of this
article is to exhibit large classes of contractions on
noncommutative $L^p$-spaces which satisfy inequality
(\ref{LpMatsaev}) for any complex polynomial $P$. The next theorem
gathers some of our main results.

\begin{thm}
Suppose $1 <p <\infty $. The following maps satisfy
(\ref{LpMatsaev}) for any complex polynomial $P$.
\begin{enumerate}
  \item A Schur multiplier $S^p \xra{M_A} S^p$ induced by a contractive Schur multiplier $B\big(\ell^2\big) \xra{M_A} B\big(\ell^2\big)$ associated with a real-valued matrix $A$.

  \item A Fourier multiplier $L^p\big(\VN(G)\big)\xra{M_t} L^p\big(\VN(G)\big)$ induced by a contractive Fourier multiplier
        $\VN(G)\xra{M_t} \VN(G)$ associated with a real valued function $G\xra{t}\mathbb{R}$, in the case where $G$
        is an amenable discrete group $G$.

  \item A Fourier multiplier $L^p\big(\VN(\mathbb{F}_n)\big)
        \xra{M_t} L^p\big(\VN(\mathbb{F}_n)\big)$ induced by a unital completely positive Fourier multiplier $\VN(\mathbb{F}_n)\xra{M_t} \VN(\mathbb{F}_n)$
        associated with a real valued function $\mathbb{F}_n \xra{t}\mathbb{R}$, where $\mathbb{F}_n$ is the free group with $n$ generators ($1\leq n\leq \infty$).
\end{enumerate}
\end{thm}

The proof of these results will use dilation theorems that we now
state. Moreover, these theorems rely on constructions dues to \'{E}.
Ricard \cite{Ric}.

\begin{thm}
\label{th dilatation Sp} Let $B\big(\ell^2\big) \xra{M_A}
B\big(\ell^2\big)$ be a unital completely positive Schur multiplier
with a real-valued matrix $A$. Then there exists a hyperfinite von
Neumann algebra $M$ equipped with a semifinite normal faithful
trace, a unital trace preserving $*$-automorphism $M \xra{U} M$, a
unital trace preserving one-to-one normal $*$-homomorphism
$B(\ell^2) \xra{J} M$ such that
\begin{eqnarray*}
(M_A)^k=\mathbb{E}U^kJ
\end{eqnarray*}
for any integer $k\geq 0$, where $M \xra{\mathbb{E}}
B\big(\ell^2\big)$ is the canonical faithful normal trace preserving
conditional expectation associated with $J$.
\end{thm}

\begin{thm}
\label{th dilatation Fourier multipliers} Let $G$ be a discrete
group. Let $\VN(G) \xra{M_t} \VN(G)$ be a unital completely positive
Fourier multiplier associated with a real valued function $G
\xra{t}\R$. Then there exists a von Neumann algebra $M$ equipped
with a faithful finite normal trace, a  unital trace preserving
$*$-automorphism $M \xra{U} M$, a unital normal trace preserving
one-to-one $*$-homomorphism $\VN(G) \xra{J} M$ such that,
\begin{eqnarray*}
(M_t)^k=\mathbb{E}U^kJ
\end{eqnarray*}
for any integer $k\geq 0$, where $M \xra{\mathbb{E}} VN(G)$ is the
canonical faithful normal trace preserving conditional expectation
associated with $J$. Moreover, if $G$ is amenable or if
$G=\mathbb{F}_n$ ($1\leq n\leq \infty$), the von Neumann algebra $M$
has the quotient weak expectation property.
\end{thm}

Various norms on the space of complex polynomials arise from
Matsaev's conjecture and Peller's problem, and it is interesting to
try to compare them. If $1\leq p\leq \infty$, note that the space of
all diagonal matrices in $S^p$ can be identified with $\ell^p$. In
this regard, the shift operator $\ell^p\xra{S}\ell^p$ coincides with
the restriction of $S^p\xra{\sigma}S^p$ to diagonal matrices. This
readily implies that
$$
\bnorm{P(S)}_{\ell^p \xra{}\ell^p} \leq  \bnorm{P(\sigma)}_{S^p
\xra{}S^p} \leq  \bnorm{P(\sigma)\ot Id_{S^p}}_{S^p(S^p)
\xra{}S^p(S^p)}
$$
for any complex polynomial $P$. We will show the following result,
which disproves a conjecture due to Peller \cite[Conjecture
2]{Pel2}.

\begin{thm}
\label{th exists a polynomial} Suppose $1 < p < \infty$, $p\not=2$.
Then there exists a complex polynomial $P$ such that
$$
 \bnorm{P(S)}_{\ell^p \xra{}\ell^p} <
\bnorm{P(\sigma)\ot Id_{S^p}}_{S^p(S^p) \xra{}S^p(S^p)}.
$$
\end{thm}
To complete this investigation, we will also show that
\begin{equation}\label{Equality}
\bnorm{P(\sigma)}_{S^p \xra{}S^p} = \bnorm{P(\sigma)\ot
Id_{S^p}}_{S^p(S^p) \xra{}S^p(S^p)} = \bnorm{P(S)\ot
Id_{S^p}}_{\ell^p(S^p) \xra{}\ell^p(S^p)}
\end{equation}
for any $P$ (the first of these equalities being due to \'{E}.
Ricard).

The paper is organized as follows. In \S2, we fix some notations, we
give some background on the key notion of completely bounded maps on
noncommutative $L^p$-spaces, we prove the second equality of
(\ref{Equality}) and we give some preliminary results. In \S3, we
show that some Fourier multipliers on $L^p(\R)$ and $\ell^p_\Z$ are
bounded but not completely bounded and we prove Theorem \ref{th
exists a polynomial} and the first equality of (\ref{Equality}). \S4
is devoted to classes of contractions which satisfy noncommutative
Matsaev's inequality (\ref{LpMatsaev}) for any complex polynomial
$P$. In particular we prove Theorems \ref{th dilatation Sp} and
\ref{th dilatation Fourier multipliers}. In \S5, we consider a
natural analog of Question \ref{quest matsaev semifinite} for
$C_0$-semigroups of contractions. Finally in \S6, we exhibit some
polynomials $P$ which always satisfy (\ref{LpMatsaev}) for any
contraction $T$.

%%%%%%%%%%%%%%%%%%%%%%%%%%%%%%%%%%%%%%%%%%%%%%%%%%%%%%%%%%%%%%%%
%%%%%%%%%%%%%%%%%%%%%%%%%%%%%%%%%%%%%%%%%%%%%%%%%%%%%%%%%%%%%%%%

\section{Preliminaries}

%%%%%%%%%%%%%%%%%%%%%%%%%%%%%%%%%%%%%%%%%%%%%%%%%%%%%%%%%%%%%%%%
%%%%%%%%%%%%%%%%%%%%%%%%%%%%%%%%%%%%%%%%%%%%%%%%%%%%%%%%%%%%%%%%

Let us recall some basic notations. Let $\T=\big\{z \in \C\ |\ |z|=1
\big\}$ and $\delta_{i,j}$ the symbol of Kronecker.

If $I$ is an index set and if $E$ is a vector space, we write
$\Mat_I$ for the space of the $I \times I$ matrices with entries in
$\C$ and $\Mat_I(E)$ for the space of the $I \times I$ matrices with
entries in $E$. If $K$ is another index set, we have an isomorphism
$\Mat_I(\Mat_K)=\Mat_{I \times K}$.

Let $M$ be a von Neumann algebra equipped with a semifinite normal
faithful trace $\tau$. For $1 \leq p<\infty$ the noncommutative
$L^p$-space $L^p(M)$ is defined as follows. If $S^+$ is the set of
all positive $x \in M$ such that $\tau(x)<\infty$ and $S$ is its
linear span, then $L^p(M)$ is the completion of $S$ with respect to
the norm $\norm{x}_{L^p(M)}=\tau \big(|x|^p\big)^{\frac{1}{p}}$. One
sets $L^{\infty}(M)=M$. We refer to \cite{PX}, and the references
therein, for more information on these spaces.

Let $1 \leq p< \infty$. If $I$ is an index set and if we equip the
space $B\big(\ell^2_I\big)$ with the operator norm and the canonical
trace $\tr$, the space $L^p\big(B(\ell^2_I)\big)$ identifies to the
Schatten-von Neumann class $S^p_{I}$.  The space $S^p_{I}$ is the
space of those compact operators $x$ from $\ell^2_I$ into $\ell^2_I$
such that
$\norm{x}_{S^p_{I}}=\big(\tr(x^*x)^{\frac{p}{2}}\big)^{\frac{1}{p}}<\infty$.
The space $S^\infty_{I}$ of compact operators from $\ell^2_I$ into
$\ell^2_I$ is equipped with the operator norm. For $I=\mathbb{N}$,
we simplify the notations, we let $S^p$ for $S^p_{\mathbb{N}}$.
Elements of $S^p_{I}$ are regarded as matrices $A=[a_{ij}]_{i,j \in
I}$ of $\Mat_{I}$. The space $S^p_{I}(S^p_{K})$ is the space of
those compact operators $x$ from $\ell^2_I \ot_{2} \ell^2_K$ into
$\ell^2_I \ot_{2} \ell^2_K$ such that $\norm{x}_{S^p_{I}(S^p_{K})} =
\big((\tr \ot \tr)(x^*x)^{\frac{p}{2}}\big)^{\frac{1}{p}} < \infty$.
Elements of $S^p_{I}(S^p_{K})$ are regarded as matrices of
$\Mat_{I}(\Mat_{K})$.

Let $M$ be a von Neumann algebra equipped with a semifinite normal
faithful trace $\tau$. If the von Neumann algebra
$B\big(\ell^2_I\big)\ovl{\ot}M$ is equipped with the semifinite
normal faithful trace $\tr\ot \tau$, the space
$L^p\big(B(\ell^2_I)\ovl{\ot}M\big)$ identifies to a space
$S^p_I\big(L^p(M)\big)$ of matrices of $\Mat_{I}\big(L^p(M)\big)$.
Moreover, under this identification, the algebraic tensor product
$S^p_I \ot L^p(M)$ is dense in $S^p_I\big(L^p(M)\big)$.

Let $N$ be another von Neumann algebra equipped with a semifinite
normal faithful trace. If $1\leq p \leq \infty$, we say that a
linear map $L^p(M)\xra{T}L^p(N)$ is completely bounded if $Id_{S^p}
\ot T$ extends to a bounded operator
$S^p\big(L^p(M)\big)\xra{Id_{S^p} \ot T}S^p\big(L^p(N)\big)$. In
this case, the completely bounded norm
$\norm{T}_{cb,L^p(M)\xra{}L^p(N)}$ is defined by
\begin{equation}\label{defnormecb}
\norm{T}_{cb,L^p(M)\xra{}L^p(N)}=\bnorm{Id_{S^p} \ot
T}_{S^p(L^p(M))\xra{}S^p(L^p(N))}.
\end{equation}
If $\Omega$ is a measure space, the space $S^p\big(L^p(\Omega)\big)$
is isometric to the $L^p$-space $L^p(\Omega,S^p)$ of $S^p$-valued
functions in Bochner's sense. Thus, if
$L^p(\Omega)\xra{T}L^p(\Omega)$ is a linear map, we have
\begin{equation}\label{normecbcommutatif}
\norm{T}_{cb,L^p(\Omega)\xra{}L^p(\Omega)}=\bnorm{T \ot
Id_{S^p}}_{L^p(\Omega,S^p)\xra{}L^p(\Omega,S^p)}.
\end{equation}
The notion of completely bounded map and the completely bounded norm
defined in (\ref{defnormecb}) are the same that these defined in
operator space theory, see \cite{ER1}, \cite{Pis3} and \cite{Pis4}.

Now, we let:
\begin{defi}
Let $M$ be a von Neumann algebra equipped with a faithful semifinite
normal trace and $1\leq p \leq \infty$. Let $L^p(M)\xra{T}L^p(M)$ be
a contraction. We say that $T$ satisfies the noncommutative
Matsaev's property if (\ref{LpMatsaev}) holds for any complex
polynomial $P$.
\end{defi}

We denote by $\ell^{p} \xra{S} \ell^{p}$ the right shift on
$\ell^{p}$. We use the same notation for the right shift on
$\ell^p_{\Z}$. We denote by $S_{-}$ the left shift on $\ell^{p}$
defined by $S_{-}(a_0,a_1,a_2,\ldots)=(a_1,a_2,a_3,\ldots)$. Suppose
$1\leq p \leq\infty$. Let $X$ be a Banach space. For any complex
polynomial $P$, we define $\norm{P}_{p,X}$ by
\begin{equation*}\label{defnormPX}
\norm{P}_{p,X}=\bnorm{P(S)\ot Id_{X}}_{\ell^p(X) \xra{} \ell^p(X)}.
\end{equation*}
We let
$\norm{P}_{p}=\norm{P}_{p,\mathbb{C}}\big(=\norm{P(S)}_{\ell^p
\xra{} \ell^p}\big)$. If $1\leq p<\infty$, it is easy to see that,
for any complex polynomial $P$, we have
\begin{equation}\label{PS equalities}
\norm{P}_{p,X}=\bnorm{P(S)\ot Id_{X}}_{\ell^p_{\Z}(X)\xra{}
\ell^p_{\Z}(X)}=\bnorm{P(S_{-}) \ot Id_X}_{\ell^p(X)\xra{}
\ell^p(X)}.
\end{equation}
Moreover, for all $1 \leq p < \infty$, by (\ref{normecbcommutatif}),
we have
\begin{eqnarray}\label{equality Z et N cb}
% \nonumber to remove numbering (before each equation)
 \norm{P}_{p,S^p}  &=&\bnorm{P(S)}_{cb,\ell^p_{\Z}\xra{}\ell^p_{\Z}}.
\end{eqnarray}
Note that, if $1\leq p\leq \infty$, we have
$\norm{P}_{p,S^p}=\norm{P}_{p^*,S_{p^*}}$. Moreover, if $1\leq p\leq
q\leq 2$, we have $\norm{P}_{q,S_q} \leq \norm{P}_{p,S^p}$ by
interpolation. We define the linear map $S^p_\Z \xra{\Theta}S^p_\Z$
as the shift "from NW to SE" which maps any matrix
$$
\left[
  \begin{array}{ccccc}
\cdots   &  \cdots  & \cdots & \cdots & \cdots\\
\cdots   &  a_{0,0} & a_{0,1} & a_{0,2} & \cdots\\
\cdots   &  a_{1,0}& a_{1,1} & a_{1,2} & \cdots \\
\cdots   &  a_{2,0} & a_{2,1}& a_{2,2} & \cdots\\
\cdots   &  \cdots & \cdots & \cdots & \cdots \\
  \end{array}
\right] \ \ \text{to} \ \ \left[
  \begin{array}{ccccc}
\cdots   &  \cdots  & \cdots & \cdots & \cdots\\
\cdots   &  a_{-1,-1} & a_{-1,0} & a_{-1,1} & \cdots\\
\cdots   &  a_{0,-1}& a_{0,0} & a_{0,1} & \cdots \\
\cdots   &  a_{1,-1} & a_{1,0}& a_{1,1} & \cdots\\
\cdots   &  \cdots & \cdots & \cdots & \cdots \\
  \end{array}
\right].
$$
If $1\leq p <\infty$, it is not difficult to see that for any
complex polynomial $P$ we have
\begin{equation}
\label{egalité theta sigma} \bnorm{P(\Theta)}_{S^p_{\Z}
\xra{}S^p_{\Z}} = \bnorm{P(\sigma)}_{S^p \xra{}S^p}
 \ \ \ \ \ \text{and} \ \ \ \ \
\bnorm{P(\Theta)}_{cb,S^p_{\Z} \xra{}S^p_{\Z}} =
\bnorm{P(\sigma)}_{cb,S^p \xra{}S^p}.
\end{equation}
Moreover, it is easy to see that, for all $A \in S^p_{\Z}$, we have
the equality $\Theta(A)=SAS^{-1}$ where we consider $A$ and
$\Theta(A)$ as operators on $\ell^2_\mathbb{Z}$.

We will use the following theorem inspired by a well-known technique
of Kitover.
\begin{thm}
\label{Kitover} Suppose $1 \leq p\leq \infty$. Let $X$ be a Banach
space and $X \xra{T} X$ an isometry (not necessarily  onto). For any
complex polynomial $P$, we have the inequality
$$
\bnorm{P(T)}_{X \xra{} X} \leq  \norm{P}_{p,X}.
$$
\end{thm}

\begin{preuve}
It suffices to consider the case $1< p < \infty$. Let $0<r<1$. Since
$T$ is an isometry we have
$$
\sum_{j=0}^{+\infty}\bnorm{r^jT^j(x)}_{X}^p =
\sum_{j=0}^{+\infty}r^{jp}\bnorm{T^j(x)}_{X}^p =
\norm{x}_{X}^p\Bigg(\sum_{j=0}^{+\infty}(r^{p})^j\Bigg) <+\infty.
$$
We let $\displaystyle
C_r=\Bigg(\sum_{j=0}^{+\infty}r^{jp}\Bigg)^{\frac{1}{p}}$. Now we
define the operator
$$
\begin{array}{cccc}
 W_r:   &  X   &  \longrightarrow   &  \ell^p(X)  \\
    &  x   &  \longmapsto       &  \frac{1}{C_r}\big(x,rT(x),r^2T^2(x),\ldots,r^jT^j(x),\ldots\big)  \\
\end{array}
$$
%$$
%\bnorm{W_r(x)}_{l^p(X)}
%=\frac{1}{C_r}\Bigg(\sum_{j=0}^{+\infty}r^{jp}\norm{T^j(x)}_{X}^p\Bigg)^{\frac{1}{p}}
%=\norm{x}_{X}.
%$$
which is an isometry. If $n$ is a positive integer and if $x\in X$
we have
$$
 W_r\big((rT)^nx\big)  = \frac{1}{C_r}\big(r^nT^n x,r^{n+1}T^{n+1}x,\ldots\big)   = (S_{-}\ot
 Id_X)^n\big(W_r(x)\big).
$$
We deduce that for any complex polynomial $P$ we have $ W_r
P(rT)=P(S_{-}\ot Id_X)W_r$. Now, if $x\in X$, we have
\begin{eqnarray*}
% \nonumber to remove numbering (before each equation)
 \bnorm{P(rT)x}_{X}
   &=& \bnorm{W_r\big(P(rT)x\big)}_{\ell^p(X)}  \\
   &=&  \bnorm{P(S_{-}\ot Id_X)W_r(x)}_{\ell^p(X)} \\
%   &\leq& \bnorm{P(\mathcal{S}\ot Id_X)}_{l^p(X) \xra{} l^p(X)}\ \bnorm{W_r(x)}_{l^p(X)}\\
   &\leq& \bnorm{P(S_{-})\ot Id_X}_{\ell^p(X) \xra{} \ell^p(X)}\ \norm{x}_{X}\\
   &=& \norm{P}_{p,X} \norm{x}_{X} \hspace{1cm} \text{by (\ref{PS equalities})}.
\end{eqnarray*}
Consequently, letting $r$ to $1$, we obtain finally that
$\bnorm{P(T)}_{X \xra{} X} \leq\norm{P}_{p,X}$.
\end{preuve}

\begin{cor}
\label{sigmacb=Scb} Suppose $1 \leq p \leq \infty$. Let $P$ be a
complex polynomial. We have
$$
\bnorm{P(\sigma)\ot Id_{S^p}}_{S^p(S^p) \xra{}S^p(S^p)} =
\bnorm{P}_{p,S^p}.
$$
\end{cor}

\begin{preuve}
With the diagonal embedding of $\ell^p$ in $S^p$, we see that for
any complex polynomial $P$ we have
$$
\bnorm{P}_{p,S^p} \leq\bnorm{P(\sigma)\ot Id_{S^p}}_{S^p(S^p)
\xra{}S^p(S^p)}.
$$
Now the map $S^p(S^p)\xra{\sigma \ot Id_{S^p}} S^p(S^p)$ is an
isometry. Hence, by the above theorem, we deduce that for every
complex polynomial $P$ we have
$$
\bnorm{P(\sigma)\ot Id_{S^p}}_{S^p(S^p) \xra{}S^p(S^p)}=
\bnorm{P(\sigma\ot Id_{S^p})}_{S^p(S^p) \xra{}
S^p(S^p)}\leq\norm{P}_{p,S^p(S^p)}=\norm{P}_{p,S^p}.
$$
\end{preuve}

Let $M$ be a von Neumann algebra. Let us recall that $M$ has QWEP
means that $M$ is the quotient of a $C^*$-algebra having the weak
expectation property (WEP) of C. Lance (see \cite{Oza} for more
information on these notions). It is unknown whether every von
Neumann algebra has this property. We will need the following
theorem which is a particular case of a result of \cite{Jun}.

\begin{thm}
\label{thJunge} Let $M$ be a von Neumann algebra with $\rm{QWEP}$
equipped with a faithful semifinite normal trace. Suppose $1< p <
\infty$. Let $\Omega$ be a measure space. Suppose that
$L^p(\Omega)\xra{T} L^p(\Omega)$ is a completely bounded map. Then
$T \ot Id_{L^p(M)}$ extends to a bounded operator and we have
$$
\bnorm{T \ot
Id_{L^p(M)}}_{L^p(\Omega,L^p(M))\xra{}L^p(\Omega,L^p(M))}\leq
\norm{T}_{cb,L^p(\Omega)\xra{}L^p(\Omega)}.
$$
\end{thm}
In the case where $M$ is a hyperfinite von Neumann algebra, the
statement of this theorem is easy to prove $\big($use
\cite[(3.1)]{Pis3} and \cite[(3.6)]{Pis3}$\big)$. With this theorem,
we deduce the following proposition.
\begin{prop}
\label{PSQWEP} Suppose $1 < p < \infty$. Let $M$ be a von Neumann
algebra with $\rm{QWEP}$ equipped with a faithful semifinite normal
trace. For all complex polynomial $P$ we have
$$
\norm{P}_{p,L^p(M)}\leq \norm{P}_{p,S^p}.
$$
\end{prop}
With this proposition, we can prove the following corollary.
\begin{cor}
\label{dilationQWEP} Let $M$ be a von Neumann algebra equipped with
a faithful semifinite normal trace and $1<p<\infty$. Let
$L^p(M)\xra{T}L^p(M)$ be a contraction. Suppose that there exists a
von Neumann algebra $M$ with QWEP equipped with a faithful
semifinite normal trace, an isometric embedding
$L^p(M)\xra{J}L^p(N)$, an isometry $L^p(N)\xra{U}L^p(N)$ and a
contractive projection $L^p(N)\xra{Q}L^p(M)$ such that,
\begin{eqnarray*}
T^k=QU^kJ
\end{eqnarray*}
for any integer $k\geq 0$. Then the contraction $T$ has the
noncommutative Matsaev's property.
\end{cor}

\begin{preuve}
For any complex polynomial $P$, we have
\begin{eqnarray*}
% \nonumber to remove numbering (before each equation)
\bnorm{P(T)}_{L^p(M) \xra{} L^p(M) }
   &=&\bnorm{QP(U)J}_{L^p(M) \xra{} L^p(M) }  \\
   &\leq& \bnorm{P(U)}_{L^p(N)\xra{}L^p(N)}.
\end{eqnarray*}
By using Theorem \ref{Kitover}, we obtain the inequality
$$
\bnorm{P(T)}_{L^p(M)\xra{} L^p(M) } \leq \norm{P}_{p, L^p(N)}.
$$
Now, the von Neumann algebra $N$ is QWEP. Then, by Proposition
\ref{PSQWEP}, we obtain finally that
$$
\bnorm{P(T)}_{L^p(M) \xra{} L^p(M) } \leq \norm{P}_{p,S^p}.
$$
\end{preuve}

Theorem \ref{thJunge}, Proposition \ref{PSQWEP} and Corollary
\ref{dilationQWEP} hold true more generally for noncommutative
$L^p$-spaces of a von Neumann algebra equipped with a distinguished
normal faithful state $M\xra{\varphi}\C$, constructed by Haagerup.
See \cite{PX} and the references therein  for more informations on
these spaces.

We refer to \cite{ALM}, \cite{AkS}, \cite{JLM} and \cite{Pel1} for
information on dilations on $L^p$-spaces (commutative and
noncommutative).

%%%%%%%%%%%%%%%%%%%%%%%%%%%%%%%%%%%%%%%%%%%%%%%%%%%%%%%%%%%%%%%%%%%%%%%%%%%%%%%%%%%%%%%%%%%%%%%%%%%%%%%%%%%%%%%%%%%%%%%%%%%%%%%%
%%%%%%%%%%%%%%%%%%%%%%%%%%%%%%%%%%%%%%%%%%%%%%%%%%%%%%%%%%%%%%%%%%%%%%%%%%%%%%%%%%%%%%%%%%%%%%%%%%%%%%%%%%%%%%%%%%%%%%%%%%%%%%%%
\section{Comparison between the commutative and noncommutative cases}

%%%%%%%%%%%%%%%%%%%%%%%%%%%%%%%%%%%%%%%%%%%%%%%%%%%%%%%%%%%%%%%%%%%%%%%%%%%%%%%%%%%%%%%%%%%%%%%%%%%%%%%%%%%%%%%%%%%%%%%%%%%%%%%%
%%%%%%%%%%%%%%%%%%%%%%%%%%%%%%%%%%%%%%%%%%%%%%%%%%%%%%%%%%%%%%%%%%%%%%%%%%%%%%%%%%%%%%%%%%%%%%%%%%%%%%%%%%%%%%%%%%%%%%%%%%%%%%%%
Suppose $1< p < \infty$. Let $G$ be a locally compact abelian group
with dual group $\widehat{G}$. An operator $L^p(G) \xra{T} L^p(G)$
is a Fourier multiplier if there exists a function $\psi \in
L^\infty\big(\widehat{G}\big)$ such that for any $f \in L^p(G)\cap
L^2(G)$ we have $\mathcal{F}\big(T(f)\big)= \psi\mathcal{F}(f)$
where $\mathcal{F}$ denotes the Fourier transform. In this case, we
let $T=M_{\psi}$. G. Pisier showed that, if $G$ is a compact group
and $1 < p<\infty$, $p\not=2$, there exists a bounded Fourier
multiplier $L^p(G)\xra{T}L^p(G)$ which is not completely bounded
(see \cite[Proposition 8.1.3]{Pis3}. We will show this result is
also true for the groups $\R$ and $\Z$ and we will prove Theorem
\ref{th exists a polynomial}.

If $b\in L^1(G)$, we define the convolution operator $C_b$ by
$$
\begin{array}{cccc}
  C_b:  &  L^p(G)   &  \longrightarrow   & L^p(G)   \\
        &   f  &  \longmapsto       &  b*f.  \\
\end{array}
$$
This operator is a completely bounded Fourier multiplier. We observe
that, if $P=\sum_{k=0}^{n}a_kz^k$ is a complex polynomial, the
operator $\ell^p_{\Z} \xra{P(S)}\ell^p_{\Z}$ is the operator
$\ell^p_{\Z} \xra{C_{\tilde{a}}} \ell^p_{\Z}$ where $\tilde{a}$ is
the sequence defined by $\tilde{a}_k=a_k$ if $0 \leq k\leq n$ and
$\tilde{a}_k=0$ otherwise.

We  will use the following approximation result \cite[Theorem
5.6.1]{Lar}.

\begin{thm}\label{th Larsen}
Suppose $1\leq p<\infty$. Let $G$ be a locally compact abelian
group. Let $L^p(G) \xra{T} L^p(G)$  be a bounded Fourier multiplier.
Then there exists a net of continuous functions $(b_l)_{i\in L}$
with compact support such that
$$
\bnorm{C_{b_{l}}}_{L^p(G) \xra{}L^p(G)}\leq \norm{T}_{L^p(G)
\xra{}L^p(G)} \ \ \ \text{and} \ \ \ C_{b_{l}} \xra[l]{so} T
$$
(convergence for the strong operator topology).
\end{thm}
Moreover, we need the following vectorial extension of
\cite[Proposition 3.3]{DeL}. One can prove this theorem as
\cite[Theorem 3.4]{CoW}.
\begin{thm}\label{deleuw}
Suppose $1 < p < \infty$. Let $\psi$ be a continuous function on
$\R$ which defines a completely bounded Fourier multiplier $M_\psi$
on $L^{p}(\R)$. Then the restriction $\psi|\Z$ of the function
$\psi$ to $\Z$ defines a completely bounded Fourier multiplier
$M_{\psi|\Z}$ on $L^{p}(\T)$.%and we
%have
%$$
%\bnorm{M_{(\psi(k))_{k \in \Z}}}_{cb,L^p(\T) \xra{} L^p(\T)} \leq
%\bnorm{M_\psi}_{cb,L^p(\R) \xra{} L^p(\R)}.
%$$
\end{thm}
We will use the next result of Jodeit \cite[Theorem 3.5]{Jod}. We
introduce the function $\Lambda\colon \mathbb{R} \xra{} \mathbb{R}$
defined by
$$
\Lambda(x)=\bigg\{\begin{array}{cl}
        1-|x|  &  {\rm if}\quad  x \in [-1,1]\\
         0     &  {\rm if}\quad  |x|>1.
       \end{array}
$$
\begin{thm}
\label{Jodeit} Suppose $1<p<\infty$.  Let $\varphi$ be a complex
function defined on $\Z$ such that $M_{\varphi}$ is a bounded
Fourier multiplier on $L^p(\T)$. Then the complex function $
\mathbb{R} \xra{\psi} \mathbb{C}$ defined on $\mathbb{R}$ by
\begin{equation}\label{def psi}
\psi(x)=\sum_{k\in\mathbb{Z}} \varphi(k)\Lambda(x-k),\qquad x\in
\mathbb{R},
\end{equation}
defines a bounded Fourier multiplier
$L^p(\mathbb{R})\xra{M_{\psi}}L^p(\mathbb{R})$.
\end{thm}
Now, we are ready to prove the following theorem.
\begin{thm}\label{mult bounded not completely bounded on R}
Suppose $1<p<\infty$, $p\not=2$. Then there exists a bounded Fourier
multiplier $L^p(\mathbb{R})\xra{M_{\psi}}L^p(\mathbb{R})$ which is
not completely bounded.
\end{thm}

\begin{preuve}
By \cite[Proposition 8.1.3]{Pis3}, there exists a bounded Fourier
multiplier $L^p(\mathbb{T}) \xra{M_{\varphi}} L^p(\mathbb{T})$ which
is not completely bounded. Now, we define the function $\psi$ on
$\mathbb{R}$ by (\ref{def psi}). By Theorem \ref{Jodeit}, the
function $\mathbb{R}\xra{\psi}\mathbb{C}$ defines a bounded Fourier
multiplier $L^p(\mathbb{R}) \xra{M_{\psi}} L^p(\mathbb{R})$. Now,
suppose that $M_{\psi}$ is completely bounded. Since the function
$\mathbb{R}\xra{\psi}\mathbb{C}$ is continuous, by Theorem
\ref{deleuw}, we deduce that the restriction $\psi|\Z$ defines a
completely bounded Fourier multiplier $M_{\psi|\Z}$ on
$L^p(\mathbb{T})$. Moreover, we observe that, for all $k \in
\mathbb{Z}$, we have
$$
\psi(k)=\varphi(k).
$$
Then we deduce that the Fourier multiplier $L^p(\mathbb{T})
\xra{M_{\varphi}} L^p(\mathbb{T})$ is completely bounded. We obtain
a contradiction. Consequently, the bounded Fourier multiplier
$L^p(\mathbb{R}) \xra{M_{\psi}} L^p(\mathbb{R})$ is not completely
bounded.
\end{preuve}

The proof of the next theorem is inspired by \cite[page 25]{CoW}.

\begin{thm}
\label{mult bounded not completely bounde} Suppose $1<p<\infty$,
$p\not=2$. Then
\begin{enumerate}
  \item There exists a bounded Fourier multiplier $\ell^p_\mathbb{Z}\xra{T}\ell^p_\mathbb{Z}$ which is not completely bounded.
  \item There exists a complex polynomial $P$ such that $\norm{P}_{p}<\norm{P}_{p,S^p}$.
\end{enumerate}
\end{thm}

\begin{preuve}
By Theorem \ref{mult bounded not completely bounded on R}, there
exists a bounded Fourier multiplier
$L^p(\mathbb{R})\xra{M_{\psi}}L^p(\mathbb{R})$ which is not
completely bounded. We can suppose that $M_{\psi}$ satisfies
$\bnorm{M_{\psi}}_{L^p(\mathbb{R}) \xra{}L^p(\mathbb{R})}=1$. By
Theorem \ref{th Larsen}, there exists a net of continuous functions
$(b_l)_{l \in L}$ with compact support such that
$$
\bnorm{C_{b_{l}}}_{L^p(\mathbb{R}) \xra{} L^p(\mathbb{R})}\leq 1 \ \
\ \text{and} \ \ \ C_{b_{l}} \xra[l]{so} M_{\psi}.
$$
Let $c>1$. There exists an element $y=\sum_{k=1}^n f_k \ot x_k\in
L^p(\mathbb{R}) \ot S^p$ with $\norm{y}_{L^p(\mathbb{R},S^p)} \leq
1$ such that $\bnorm{(M_{\psi}\ot
Id_{S^p})(y)}_{L^p(\mathbb{R},S^p)}\geq 3c$. Then, it is not
difficult to see that there exists $l\in L$ such that
$\bnorm{(C_{b_{l}}\ot Id_{S^p})(y)}_{L^p(\mathbb{R},S^p)}\geq 2c$.
We deduce that there exists a continuous function $b\colon
\mathbb{R} \xra{}\mathbb{C}$ with compact support such that
$\norm{C_b}_{L^p(\mathbb{R}) \xra{} L^p(\mathbb{R})} \leq 1$ and
$\bnorm{C_b}_{cb,L^p(\mathbb{R}) \xra{} L^p(\mathbb{R})} \geq 2c$.
%Let $c>1$. Now, suppose that, for any $l\in L$, the completely
%bounded norm of the Fourier multiplier $C_{b_{l}}$ satisfies
%$\bnorm{C_{b_{l}}}_{cb,L^p(\mathbb{R}) \xra{}L^p(\mathbb{R})}\leq
%2c$. Then, for any matrix $A=[a_{ij}]_{1\leq i,j\leq n}\in
%M_n\big(L^p(\mathbb{R})\big)$, we have
%\begin{eqnarray*}
%% \nonumber to remove numbering (before each equation)
%\Bnorm{\big[M_{\psi}(a_{ij})\big]_{1\leq i,j\leq
%n}}_{M_n(L^p(\mathbb{R}))}
%   &=& \lim_{l} \Bnorm{\big[C_{b_{l}}(a_{ij})\big]_{1\leq i,j\leq n}}_{M_n(L^p(\mathbb{R}))} \\
%   &\leq& 2c\norm{A}_{M_n(L^p(\mathbb{R}))}
%\end{eqnarray*}
%(here we use the canonical operator space structure on $L^p(G)$, see
%\cite{Pis3}).
%Therefore we deduce that the linear map
%$L^p(\mathbb{R})\xra{M_{\psi}}L^p(\mathbb{R})$ is completely
%bounded. We obtain a contradiction.
Thus there exists a continuous function $\mathbb{R}
\xra{b}\mathbb{C}$ with compact support such that
$$
\norm{C_b}_{L^p(\mathbb{R}) \xra{} L^p(\mathbb{R})} \leq 1 \ \ \ \ \
\text{and} \ \ \ \ \ \norm{C_b}_{cb,L^p(\mathbb{R}) \xra{}
L^p(\mathbb{R})} \geq 2c.
$$
Now, we define the sequence $\big(a_n\big)_{n\geq 1}$ of complex
sequences indexed by $\mathbb{Z}$ by, if $n\geq 1$ and $k \in
\mathbb{Z}$
$$
a_{n,k}=\int_{0}^{1}\int_{0}^{1}\frac{1}{n}b\bigg(\frac{t-s+k}{n}\bigg)dsdt.%=\int_{-\frac{1}{n}}^{\frac{1}{n}}b\bigg(\frac{k}{n}+s\bigg)\big(1-n|s|\big)ds
$$
Note that each sequence $a_{n}$ has only a finite number of non-zero
term. Let $n\geq 1$. We introduce the conditional expectation
$L^p(\mathbb{R}) \xra{\mathbb{E}_n} L^p(\mathbb{R})$ with respect to
the $\sigma$-algebra generated by the
$\Big[\frac{k}{n},\frac{k+1}{n}\Big[$, $k \in \mathbb{Z}$. For every
integer $n\geq 1$ and all $f \in L^p(\mathbb{R})$, we have
%\vspace{0.3cm}
$$
\mathbb{E}_nf=n \sum_{k \in
\mathbb{Z}}^{}\Bigg(\int_{\frac{k}{n}}^{\frac{k+1}{n}} f(t)dt\Bigg)
1_{\big[\frac{k}{n},\frac{k+1}{n}\big[}
$$
%\vspace{0.4cm} \noindent
(see \cite[page 227]{AbA}). Now, we define the linear map
$\ell^p_\mathbb{Z} \xra{J_n} \mathbb{E}_n\big(L^p(\mathbb{R})\big)$
by, if $u \in \ell^p_\mathbb{Z}$ \vspace{0.2cm}
$$
J_n(u)=n^{\frac{1}{p}}\sum_{k\in \mathbb{Z}} u_k
1_{\big[\frac{k}{n},\frac{k+1}{n}\big[}.
$$
%\vspace{0.4cm} \noindent
It is easy to check that the map $J_n$ is an isometry of
$\ell^p_\mathbb{Z}$ onto the range
$\mathbb{E}_n\big(L^p(\mathbb{R})\big)$ of $\mathbb{E}_n$. For any
$u\in \ell^p_\mathbb{Z}$, we have
\begin{align}
% \nonumber to remove numbering (before each equation)
\mathbb{E}_nC_bJ_n(u)   &= n \sum_{k \in \mathbb{Z}}^{}\Bigg(\int_{\frac{k}{n}}^{\frac{k+1}{n}}\big(C_bJ_n(u)\big)(t)dt\Bigg)1_{\big[\frac{k}{n},\frac{k+1}{n}\big[}\nonumber \\
   &= n \sum_{k \in \mathbb{Z}}^{}\Bigg(\int_{\frac{k}{n}}^{\frac{k+1}{n}}\int_{-\infty}^{+\infty}b(t-s)\big(J_n(u)\big)(s)dsdt\Bigg)1_{\big[\frac{k}{n},\frac{k+1}{n}\big[} \nonumber\\
   &= n \sum_{k \in \mathbb{Z}}^{}\Bigg(\int_{\frac{k}{n}}^{\frac{k+1}{n}}\int_{-\infty}^{+\infty}b(t-s)n^{\frac{1}{p}}\Bigg(\sum_{j \in\mathbb{Z}} u_{j} 1_{\big[\frac{j}{n},\frac{j+1}{n}\big[}(s)\Bigg)dsdt\Bigg)1_{\big[\frac{k}{n},\frac{k+1}{n}\big[} \label{summation}\\
   &= n^{1+\frac{1}{p}} \sum_{k \in \mathbb{Z}}^{}\Bigg(\sum_{j \in\mathbb{Z}}u_{j}\int_{\frac{k}{n}}^{\frac{k+1}{n}}\int_{\frac{j}{n}}^{\frac{j+1}{n}}b(t-s)dsdt\Bigg)1_{\big[\frac{k}{n},\frac{k+1}{n}\big[} \label{equality summation} \\
   &= n^{1+\frac{1}{p}} \sum_{k \in \mathbb{Z}}^{}\Bigg(\sum_{j \in\mathbb{Z}}u_{j}\int_{0}^{\frac{1}{n}}\int_{0}^{\frac{1}{n}}b\bigg(t-s+\frac{k-j}{n}\bigg)dsdt\Bigg)1_{\big[\frac{k}{n},\frac{k+1}{n}\big[} \nonumber\\
   &= n^{\frac{1}{p}} \sum_{k \in \mathbb{Z}}^{}\Bigg(\sum_{j \in\mathbb{Z}}u_{j}\int_{0}^{1}\int_{0}^{1}b\bigg(\frac{t-s+k-j}{n}\bigg)dsdt\Bigg)1_{\big[\frac{k}{n},\frac{k+1}{n}\big[}\nonumber \\
   &= n^{\frac{1}{p}} \sum_{k \in \Z}^{} \Bigg(\sum_{j\in \Z}u_{j}a_{n,k-j}\Bigg)1_{\big[\frac{k}{n},\frac{k+1}{n}\big[}\nonumber\\
   &= J_nC_{a_n}(u)\nonumber
\end{align}
%\noindent
(where the equality (\ref{equality summation}) follows from the fact
that the summation over $j\in \Z$ of (\ref{summation}) is finite).
Thus we have the following commutative diagram
$$
 \xymatrix @R=0.58cm @C=3cm{
         L^p(\mathbb{R})   \ar[r]^{C_b}     & L^p(\mathbb{R})\ar[d]^{\mathbb{E}_n}\\
         \mathbb{E}_n\big(L^p(\mathbb{R})\big)   \ar@{^{(}->}[u]               &\mathbb{E}_n\big(L^p(\mathbb{R})\big)\\
         \ell^p_\mathbb{Z}  \ar[u]^{J_n}_{\approx} \ar[r]_{C_{a_n}}              &
         \ell^p_\mathbb{Z}\ar[u]_{J_n}^{\approx}.
}
$$
Then, for any integer $n\geq 1$, since
$\norm{\mathbb{E}_n}_{L^p(\mathbb{R})\xra{}L^p(\mathbb{R})}\leq 1$,
we have the following estimate
$$
\bnorm{C_{a_n}}_{\ell^p_\mathbb{Z} \xra{}\ell^p_\mathbb{Z}} \leq
\norm{C_b}_{ L^p(\mathbb{R}) \xra{} L^p(\mathbb{R})}\leq 1.
$$
Moreover, we have $\mathbb{E}_n \ot Id_{S^p}\xra[n \to
+\infty]{so}Id_{L^p(\mathbb{R},S^p)}$ (see \cite[Theorem 1]{Cha}).
It is easy to see that
$$
(\mathbb{E}_nC_b\mathbb{E}_n) \ot Id_{S^p} \xra[n \to +\infty]{so}
C_b \ot Id_{S^p}.
$$
By the strong semicontinuity of the norm, we obtain that
$$
\norm{C_b}_{cb,L^p(\mathbb{R}) \xra{} L^p(\mathbb{R})} \leq
\liminf_{n \to \infty}
\bnorm{\mathbb{E}_nC_b\mathbb{E}_n}_{cb,L^p(\mathbb{R})
\xra{}L^p(\mathbb{R})}.
$$
%By an argument of approximation, we deduce that
%$$
%\bnorm{C_b}_{cb,L^p(\mathbb{R}) \xra{} L^p(\mathbb{R})} \leq
%\liminf_{n \to \infty} \bnorm{C_{a_n}}_{cb,\ell^p(\mathbb{Z})
%\xra{}\ell^p(\mathbb{Z})}.
%$$
Then, there exists an integer $n\geq 1$ such that
$$
\bnorm{C_{a_n}}_{\ell^p_\mathbb{Z} \xra{}\ell^p_\mathbb{Z}} \leq 1 \
\ \ \ \ \text{and} \ \ \ \ \ \bnorm{C_{a_n}}_{cb,\ell^p_\mathbb{Z}
\xra{}\ell^p_\mathbb{Z}} \geq c.
$$
Thus, we prove the second assertion by shifting the obtained
multiplier. Finally, we show the first assertion by the closed graph
theorem, (\ref{PS equalities}) for $X=\mathbb{C}$ and (\ref{equality
Z et N cb}).
\end{preuve}

The paper \cite{Arh} is a continuation of these investigations. The
author proves that if $G$ is an arbitrary infinite locally compact
abelian group, $1<p<\infty$ and $p\not=2$ then there exists a
bounded Fourier multiplier on $L^p(G)$ which is not completely
bounded.

In the light of Corollary \ref{sigmacb=Scb} and Theorem \ref{mult
bounded not completely bounde}, it is natural to compare
$\bnorm{P(\sigma)}_{cb,S^p \xra{}S^p}$ and $\bnorm{P(\sigma)}_{S^p
\xra{}S^p}$. We finish the section by proving that these quantities
are identical. It is a result due to \'{E}. Ricard. In order to
prove it, we need the following notion of Schur multiplier. We equip
$\T$ with its normalized Haar measure. We denote by
$S^p\big(L^2(\T)\big)$ the Schatten-von Neumann class associated
with $B\big(L^2(\T)\big)$. If $f\in L_{2}(\T\times \T)$, we denote
the associated Hilbert-Schmidt operator by
$$
\begin{array}{cccc}
  K_f : &  L^2(\T)   &  \longrightarrow   & L^2(\T)   \\
    &   u  &  \longmapsto       &  \int_{\T}u(z)f(z,\cdot)dz.  \\
\end{array}
$$
A Schur multiplier on $S^p\big(L^2(\T)\big)$ is a linear map
$S^p\big(L^2(\T)\big)\xra{T}S^p\big(L^2(\T)\big)$ such that there
exists a measurable function $\T \times \T \xra{\varphi}\C$ which
satisfies, for any finite rank operator of the form
$L^2(\T)\xra{K_f}L^2(\T)$, the equality $T(K_f)=K_{\varphi f}$. We
denote $T$ by $M_{\varphi}$ and we say that the function $\varphi$
is the symbol of the Schur multiplier
$S^p\big(L^2(\T)\big)\xra{M_{\varphi}}S^p\big(L^2(\T)\big)$ (see
\cite{BiS} and \cite{LaS} for more details).

We denote by $L^2(\T)\xra{\mathcal{F}}\ell^2_{\Z}$ the Fourier
transform. We define the isometry $\Psi$ by
$$
\begin{array}{cccc}
 \Psi:   &  S^p\big(L^2(\T)\big)   &  \longrightarrow   &  S^p_\Z \\
         &   T  &  \longmapsto       & \mathcal{F}T\mathcal{F}^{-1}.   \\
\end{array}
$$
Now, we can show the following proposition.
\begin{prop}
Suppose $1 \leq p \leq \infty$. For any complex polynomial $P$, we
have
$$
\bnorm{P(\sigma)}_{S^p \xra{}S^p}=\bnorm{P(\sigma)}_{cb,S^p
\xra{}S^p}\Big(=\bnorm{P(\sigma)\ot Id_{S^p}}_{S^p(S^p)
\xra{}S^p(S^p)}\Big).
$$
\end{prop}

\begin{preuve}
It suffices to consider the case $1<p<\infty$. For any $n\in \Z$ and
any finite rank operator of the form $K_f$, we have
\begin{align*}
% \nonumber to remove numbering (before each equation)
 \Big(\Theta\Psi\big(K_f\big)\Big)(e_n)
   &=  S\mathcal{F}K_f\mathcal{F}^{-1}S^{-1}(e_n)  \\
   &= S\mathcal{F}K_f(z^{n-1}) \\
   &= S\mathcal{F}\Bigg(\int_{\T} z^{n-1}f(z,\cdot)dz\Bigg) \\
   &= S\Bigg(\sum_{k\in \Z}\Bigg(\int_{\T^2} z^n\ovl{z'}^kf(z,z')dzdz'\Bigg)e_k\Bigg) \\
   &= \sum_{k\in \Z}\Bigg(\int_{\T^2} z^{n-1} \ovl{z'}^k f(z,z')dzdz'\Bigg)e_{k+1}.
\end{align*}
Now we define the function $\T \times \T \xra{\varphi}\C$ by
$\varphi(z,z')=z^{-1}z'$ where $z,z'\in \T$. Then, for any $n\in \Z$
and any finite rank operator of the form $K_f$, we have
\begin{align*}
% \nonumber to remove numbering (before each equation)
 \Big(\Psi M_\varphi\big(K_f\big)\Big)(e_n)  %&=& \Big(\Psi\big(K_{\varphi f}\big)\Big)(e_n) \\
   &= \mathcal{F}K_{\varphi f}\mathcal{F}^{-1}(e_n)\\
   &= \mathcal{F}\Bigg(\int_{\T}z^n\varphi(z,\cdot)f(z,\cdot)dz\Bigg)\\
   &= \sum_{k\in \Z}\Bigg(\int_{\T^2} z^{n} \ovl{z'}^k \varphi(z,z')f(z,z')dzdz'\Bigg)e_{k}\\
   &= \sum_{k\in \Z}\Bigg(\int_{\T^2} z^{n-1} \ovl{z'}^{k-1} f(z,z')dzdz'\Bigg)e_{k}\\
   &= \sum_{k\in \Z}\Bigg(\int_{\T^2} z^{n-1} \ovl{z'}^k f(z,z')dzdz'\Bigg)e_{k+1}.
\end{align*}

%\begin{eqnarray*}
%% \nonumber to remove numbering (before each equation)
% \Big(\Psi M_\varphi\big(K_f\big)\Big)(e_n)  %&=& \Big(\Psi\big(K_{\varphi f}\big)\Big)(e_n) \\
%   &=& \mathcal{F}K_{\varphi f}\mathcal{F}^{-1}(e_n)\\
%   &=& \mathcal{F}\Bigg(\int_{\T}z^n\varphi(z,\cdot)f(z,\cdot)dz\Bigg)\\
%   &=& \sum_{k\in \Z}\Bigg(\int_{\T^2} z^{n} \ovl{z'}^k \varphi(z,z')f(z,z')dzdz'\Bigg)e_{k}
%\end{eqnarray*}
%
%\begin{eqnarray*}
%\ \ \ \ \ \ \ \ \ \ \ \ \ \ \ \ \   &=& \sum_{k\in \Z}\Bigg(\int_{\T^2} z^{n-1} \ovl{z'}^{k-1} f(z,z')dzdz'\Bigg)e_{k}\\
%&=& \sum_{k\in \Z}\Bigg(\int_{\T^2} z^{n-1} \ovl{z'}^k
%f(z,z')dzdz'\Bigg)e_{k+1}.
%\end{eqnarray*}
Then, for any complex polynomial $P$, we have the following
commutative diagram
$$
 \xymatrix @R=1cm @C=1cm{
         S^p_\Z                                   \ar[r]^{P(\Theta)} &  S^p_\Z                         \\
       S^p\big(L^2(\T)\big) \ar[u]^{\Psi}   \ar[r]_{P(M_{\varphi})} & S^p\big(L^2(\T)\big) \ar[u]_{\Psi}. }
$$
Furthermore, for any complex polynomial $P$, we have
$P(M_{\varphi})=M_{P(\varphi)}$. Moreover, the Schur multiplier
$S^p\big(L^2(\T)\big)\xra{M_{P(\varphi)}} S^p\big(L^2(\T)\big)$ has
a continuous symbol whose the support has no isolated point. By
\cite[Theorem 1.19]{LaS}, we deduce that the norm and the completely
bounded norm of $P(M_{\varphi})$ coincide. Since $\Psi$ is a
complete isometry, we obtain the result by (\ref{egalité theta
sigma}).
\end{preuve}

%%%%%%%%%%%%%%%%%%%%%%%%%%%%%%%%%%%%%%%%%%%%%%%%%%%%%%%%%%%%%%%%
%%%%%%%%%%%%%%%%%%%%%%%%%%%%%%%%%%%%%%%%%%%%%%%%%%%%%%%%%%%%%%%%

\section{Positive results}

%%%%%%%%%%%%%%%%%%%%%%%%%%%%%%%%%%%%%%%%%%%%%%%%%%%%%%%%%%%%%%%%
%%%%%%%%%%%%%%%%%%%%%%%%%%%%%%%%%%%%%%%%%%%%%%%%%%%%%%%%%%%%%%%%
Let $M$ and $N$ be von Neumann algebras equipped with faithful
semifinite normal traces $\tau_M$ and $\tau_N$. Let $M\xra{T}N$ a
positive linear map. We say that $T$ is trace preserving if for all
$x\in L^1(M)\cap M_+$ we have $\tau_N\big(T(x)\big)=\tau_M(x)$. We
will use the following straightforward extension of \cite[Lemma
1.1]{JuX}.
\begin{lemma}
\label{lemma trace preserving} Let $M$ and $N$ be von Neumann
algebras equipped with faithful semifinite normal traces. Let
$M\xra{T}N$  be a trace preserving unital normal positive map.
Suppose $1\leq p < \infty$. Then $T$ induces a contraction
$L^p(M)\xra{T}L^p(N)$. Moreover, if $M\xra{T}N$ is an one-to-one
normal unital $*$-homomorphism, $T$ induces an isometry
$L^p(M)\xra{T}L^p(N)$.
\end{lemma}
Let $M$ be a von Neumann algebra equipped with faithful semifinite
normal trace $\tau$ and $N$ a von Neumann subalgebra such that the
restriction of $\tau$ is still semifinite. Then, it is well-known
that the extension $L^p(M)\xra{\mathbb{E}}L^p(N)$ of the canonical
faithful normal trace preserving conditional $M \xra{\mathbb{E}} N$
is a contractive projection.

Consider the situation where $M\xra{T}M$ is a linear map such there
exists a von Neumann algebra $N$ equipped with a faithful semifinite
normal trace, a unital trace preserving $*$-automorphism $N \xra{U}
N$, a unital normal trace preserving one-to-one $*$-homomorphism $M
\xra{J} N$ such that,
\begin{eqnarray}
\label{formule puissance k} T^k=\mathbb{E}U^kJ
\end{eqnarray}
for any integer $k\geq 0$, where $N \xra{\mathbb{E}} M$ is the
canonical faithful normal trace preserving conditional expectation
associated with $J$. Then, for all $ 1\leq p < \infty$, the maps $N
\xra{U} N$ and $M \xra{J} N$ extend to isometries $L_{p}(N) \xra{U}
L_{p}(N)$ and $L^p(M) \xra{J} L^p(N)$ and the map $N
\xra{\mathbb{E}} M$ extends to a contractive projection $L^p(N)
\xra{\mathbb{E}} L^p(M)$ such that (\ref{formule puissance k}) is
also true for the induced map $L^p(M) \xra{T} L^p(M)$.

In order to prove Theorems \ref{th dilatation Sp} and \ref{th
dilatation Fourier multipliers}, we need to use fermion algebras.
Since we will study maps between $q$-deformed algebras, we recall
directly several facts about these more general algebras in the
context of \cite{BKS}. We denote by $S_n$ the symmetric group. If
$\sigma$ is a permutation of $S_n$ we denote by $|\sigma|$ the
number $\card\big\{ (i,j)\ |\ 1\leq i,j \leq n, \sigma(i)>\sigma(j)
\big\}$ of inversions of $\sigma$. Let $H$ be a real Hilbert space
with complexification $H_{\C}$. If $-1 \leq q < 1 $ the $q$-Fock
space over $H$ is
$$
\mathcal{F}_{q}(H)=\C \Omega \oplus \bigoplus_{n \geq 1}
H_{\C}^{\ot_{n}}
$$
where $\Omega$ is a unit vector, called the vacuum and where the
scalar product on $H_{\C}^{\ot_{n}}$ is given by
$$
\langle h_1 \ot \dots \ot h_n ,k_1 \ot \cdots \ot k_n
\rangle_{q}=\sum_{\sigma \in S_n} q^{|\sigma|}\langle
h_1,k_{\sigma(1)}\rangle_{H_{\C}}\cdots\langle
h_n,k_{\sigma(n)}\rangle_{H_{\C}}.
$$
If $q=-1$, we must first divide out by the null space, and we obtain
the usual antisymmetric Fock space. The creation operator $l(e)$ for
$e\in H$ is given by
$$
\begin{array}{cccc}
   l(e): &   \mathcal{F}_{q}(H)  &  \longrightarrow   &  \mathcal{F}_{q}(H)  \\
         &  h_1 \ot \dots \ot h_n   &  \longmapsto       & e \ot h_1 \ot  \dots \ot h_n.   \\
\end{array}
$$
They satisfy the $q$-relation
$$
l(f)^*l(e)-ql(e)l(f)^* = \langle f,e\rangle_{H}
Id_{\mathcal{F}_{q}(H)}.
$$
We denote by $\mathcal{F}_{q}(H)\xra{\omega(e)}\mathcal{F}_{q}(H)$
the selfadjoint operator $l(e)+l(e)^*$. The $q$-von Neumann algebra
$\Gamma_q(H)$ is the von Neumann algebra generated by the operators
$\omega(e)$ where $e\in H$. It is a finite von Neumann algebra with
the trace $\tau$ defined by
$\tau(x)=\langle\Omega,x.\Omega\rangle_{\mathcal{F}_{q}(H)}$ where
$x \in \Gamma_q(H)$.

Let $H$ and $K$ be real Hilbert spaces and $H \xra{T} K$ be a
contraction with complexification $H_{\C} \xra{T_{\C}} K_{\C}$. We
define the following linear map
$$
\begin{array}{cccc}
 \mathcal{F}_q(T):   &  \mathcal{F}_{q}(H)   &  \longrightarrow   & \mathcal{F}_{q}(K)  \\
    &  h_1 \ot \dots \ot h_n   &  \longmapsto       &    T_{\C}h_1 \ot \dots \ot T_{\C}h_n. \\
\end{array}
$$
Then there exists a unique map $ \Gamma_q(H) \xra{\Gamma_q(T)}
\Gamma_q(H) $ such that for every $x \in \Gamma_q(H)$ we have
$$
\big(\Gamma_q(T)(x)\big)\Omega=\mathcal{F}_q(T)(x\Omega).
$$
This map is normal, unital, completely positive and trace
preserving. If $H \xra{T} K $ is an isometry, $\Gamma_q(T)$ is an
injective $*$-homomorphism. If $1 \leq p < \infty$, it extends to a
contraction $L^p\big(\Gamma_q(H)\big) \xra{\Gamma_q(T)}
L^p\big(\Gamma_q(K)\big)$.

We are mainly concerned with the fermion algebra $\Gamma_{-1}(H)$.
In this case, recall that if $e \in H$ has norm 1, then the operator
$\omega(e)$ satisfies $\omega(e)^2=Id_{\mathcal{F}_{-1}(H)}$.
Moreover, we need the following Wick formula, (see \cite[page
2]{Boz} and \cite[Corollary 2.1]{EfP}). In order to state this, we
denote, if $k \geq 1$ is an integer, by $\mathcal{P}_2(2k)$ the set
of 2-partitions of the set $\{1,2,\ldots,2k\}$. If $\mathcal{V}\in
\mathcal{P}_2(2k)$ we let $c(\mathcal{V})$ the number of crossings
of $\mathcal{V}$, which is given, by the number of pairs of blocks
of $\mathcal{V}$ which cross (see \cite[page 8630]{EfP} for a
precise definition). Then, if $f_1,\ldots,f_{2k}\in H$ we have
\begin{eqnarray}\label{formule de Wick}
% \nonumber to remove numbering (before each equation)
 \tau\big(\omega(f_1)\omega(f_2)\cdots\omega(f_{2k})\big)  &=&\sum_{\mathcal{V}\in \mathcal{P}_2(2k)} (-1)^{c(\mathcal{V})}
\prod_{(i,j)\in \mathcal{V}}\langle f_i,f_j\rangle_{H}.
\end{eqnarray}
In particular, for all $e,f\in H$, we have
\begin{eqnarray}\label{petite formule de Wick}
% \nonumber to remove numbering (before each equation)
  \tau\big(\omega(e)\omega(f)\big) &=&\langle e,f \rangle_{H}.
\end{eqnarray}

Let $A=[a_{ij}]_{i,j\in I}$ be a matrix of $\mathbb{M}_I$. By
definition, the Schur multiplier on $B\big(\ell^2_I\big)$ associated
with this matrix is the unbounded linear operator $M_A$ whose domain
is the space of all $B=[b_{ij}]_{i,j\in I}$ of $B\big(\ell^2_I\big)$
such that $[a_{ij}b_{ij}]_{i,j\in I}$ belongs to
$B\big(\ell^2_I\big)$, and whose action on $B=[b_{ij}]_{i,j\in I}$
is given by $M_A(B)=[a_{ij}b_{ij}]_{i,j\in I}$. For all $i,j\in I$,
the matrix $e_{ij}$ belongs to $D(M_A)$, hence $M_A$ is densely
defined for the weak* topology. Suppose $1\leq p< \infty$. If for
any $B \in S^p_{I}$, we have $B\in D(M_A)$ and the matrix $M_{A}(B)$
represents an element of $S^p_{I}$, by the closed graph theorem, the
matrix $A$ of $\Mat_{I}$ defines a bounded Schur multiplier $S^p_{I}
\xra{M_A} S^p_{I}$. We have a similar statement for bounded Schur
multipliers on $B\big(\ell^2_I\big)$.

Recall that a matrix $A$ of $\mathbb{M}_I$ defines a contractive
Schur multiplier $B\big(\ell^2_I\big) \xra{M_A} B\big(\ell^2_I\big)$
if and only if there exists an index set $K$ and norm 1 vectors $h_i
\in \ell^2_K$ and $k_j \in \ell^2_K$ such that for all $i,j\in I$ we
have $a_{i,j}=\langle h_i,k_j\rangle_{\ell^2_K}$ (see \cite{Pau}).
If all entries of $A$ are real numbers, we can take the real vector
space $\ell^2_K(\mathbb{R})$ instead of the complex vector space
$\ell^2_K$. Finally, recall that every contractive Schur multiplier
$B\big(\ell^2_I\big) \xra{M_A} B\big(\ell^2_I\big)$ is completely
contractive (see \cite{Pau}).

We say that a matrix $A$ of $\mathbb{M}_I$ induces a completely
positive Schur multiplier $B\big(\ell^2_I\big) \xra{M_A}
B\big(\ell^2_I\big)$ if and only if for any finite set $ F\subset I$
the matrix $[a_{i,j}]_{i,j \in F}$ is positive (see \cite{Pau}). An
other well-known characterization is that there exists vectors
$h_i\in \ell^2_K(\C)$ of norm 1 such that for all $i,j\in I$ we have
$a_{i,j}=\langle h_i,h_j\rangle_{\ell^2_K}$. If $A$ is a real
matrix, we can use the real vector space $\ell^2_K(\mathbb{R})$
instead of the complex vector space $\ell^2_K$.

Let $M$ be a von Neumann algebra equipped with a semifinite normal
faithful trace $\tau$. Suppose that $M\xra{T}M$ is a normal
contraction. We say that $T$ is selfadjoint if for all $x,y \in
M\cap L^1(M)$ we have
$$
\tau\big(T(x)y^*\big)=\tau\big(x(Ty)^*\big).
$$
In this case, it is easy to see that the restriction $T|M \cap
L^1(M)$ extends to a contraction $L^1(M) \xra{T} L^1(M)$. By complex
interpolation, for any $ 1\leq p\leq \infty$, we obtain a
contractive map $L^p(M) \xra{T} L^p(M)$. Moreover, the operator
$L^2(M) \xra{T} L^2(M)$ is selfadjoint. If $M\xra{T}M$ is a normal
selfadjoint complete contraction, it is easy to see that the map
$L^p(M) \xra{T} L^p(M)$ is completely contractive for all $ 1\leq
p\leq \infty$. It is easy to see that a contractive Schur multiplier
$B\big(\ell^2_I\big) \xra{M_A} B\big(\ell^2_I\big)$ associated with
a matrix $A$ of $\mathbb{M}_I$ is selfadjoint if and only if all
entries of $A$ are real.

In order to prove the next theorem, we need the following notion of
infinite tensor product of von Neumann algebras, see \cite{Tak3}.
Given a sequence $(M_n,\tau_n)_{n\in \Z}$ of von Neumann algebras
$M_n$ equipped with faithful normal finite traces $\tau_n$, then on
the infinite minimal $C^*$-tensor product of the algebras
$(M_n)_{n\in \Z}$ there is a well-defined infinite product state $
\cdots \ot \tau_{-1} \ot \tau_0 \ot \tau_1 \ot  \cdots$. The weak
operator closure of the GNS-representation of the infinite
$C^*$-tensor product of $(M_n)_{n\in \Z}$ with respect to the state
$ \cdots \ot \tau_{-1} \ot \tau_0 \ot \tau_1 \ot  \cdots$ yields a
von Neumann algebra, called the infinite tensor product of von
Neumann algebras $M_n$ with respect to the traces $\tau_n$. We will
denote this algebra by $\ovl{\bigotimes}_{n\in \Z}(M_n,\tau_n)$. The
state $ \cdots \ot \tau_{-1} \ot \tau_0 \ot \tau_1 \ot  \cdots$
extends to a faithful normal finite trace on $\ovl{\bigotimes}_{n\in
\Z}(M_n,\tau_n)$ which we still denote by $ \cdots \ot \tau_{-1} \ot
\tau_0 \ot \tau_1 \ot \cdots$.

The following theorem states that we can dilate some Schur
multipliers. The construction (and the one of Theorem \ref{th
dilatation Fourier multipliers2}) is inspired by the work of \'{E}.
Ricard \cite{Ric}.

\begin{thm}
\label{th dilatation SpI} Let $B\big(\ell^2_I\big) \xra{M_A}
B\big(\ell^2_I\big)$ be a unital completely positive Schur
multiplier associated with a real-valued matrix $A$. Then there
exists a hyperfinite von Neumann algebra $M$ equipped with a
semifinite normal faithful trace, a unital trace preserving
$*$-automorphism $M \xra{U} M$, a unital trace preserving one-to-one
normal $*$-homomorphism $B\big(\ell^2_I\big) \xra{J} M$ such that
\begin{eqnarray*}
(M_A)^k=\mathbb{E}U^kJ
\end{eqnarray*}
for any integer $k\geq 0$, where $M \xra{\mathbb{E}}
B\big(\ell^2_I\big)$ is the canonical faithful normal trace
preserving conditional expectation associated with $J$.
\end{thm}

\begin{preuve}
Since the map $B\big(\ell^2_I\big) \xra{M_A} B\big(\ell^2_I\big)$ is
completely positive we can define a positive symmetric bilinear form
$\langle \cdot,\cdot\rangle_{\ell^{2,A}}$ on the real span of the
$e_i$, where $i\in I$, by:
\begin{eqnarray}
\label{eiej=aij}
% \nonumber to remove numbering (before each equation)
  \langle e_i,e_j\rangle_{\ell^{2,A}} &=& a_{ij}.
\end{eqnarray}
We denote by $\ell^{2,A}$ the completion of the real pre-Hilbert
obtained by quotient by the corresponding kernel. For all $i$ of $I$
we still denote by $e_i$ the class of $e_i$ in $\ell^{2,A}$. Now we
define the von Neumann algebra $M$ by
$$
M=B\big(\ell^2_I\big) \otvn \Bigg(\ovl{\bigotimes_{n\in \Z}}
\big(\Gamma_{-1}(\ell^{2,A}),\tau\big)\Bigg).
$$
Since the von Neumann algebra $\Gamma_{-1}(\ell^{2,A})$ is
hyperfinite, the von Neumann algebra $M$ is also hyperfinite. We
define the element $d$ of $M$ by
$$
d=\sum_{i\in I} e_{ii} \ot \cdots \ot I \ot \omega(e_i) \ot I \ot
\cdots
$$
where $\omega(e_i)$ is in position 0. Recall that $M_A$ is unital.
Then it is not difficult to see that $d$ is a symmetry, i.e. a
selfadjoint unitary element. We equip the von Neumann algebra $M$
with the faithful semifinite normal trace $\tau_M=\tr \ot\cdots \ot
\tau \ot \tau \ot \cdots$. We denote by $M \xra{\mathbb{E}}
B\big(\ell^2_I\big)$ the canonical faithful normal trace preserving
conditional expectation of $M$ onto $B\big(\ell^2_I\big)$. We have
$$
\mathbb{E}=Id_{B(\ell^2_I)} \ot \cdots \ot \tau \ot \tau \ot\cdots.
$$
We define the canonical injective normal unital $*$-homomorphism
$$
\begin{array}{cccc}
  J:  &  B\big(\ell^2_I\big)   &  \longrightarrow   &  M  \\
    &   x  &  \longmapsto       &   x \ot \cdots \ot I \ot I \ot \cdots.  \\
\end{array}
$$
Clearly, $J$ preserves the traces. We define the right shift
%\pagebreak[4]
$$
\begin{array}{cccc}
   \mathscr{S}: &   \ovl{\bigotimes}_{n\in \Z}\big( \Gamma_{-1}(\ell^{2,A}),\tau\big)  &  \longrightarrow   & \ovl{\bigotimes}_{n\in \Z}\big(\Gamma_{-1}(\ell^{2,A}),\tau\big) \\
   &    \cdots \ot x_0 \ot x_1 \ot \cdots                       &  \longmapsto       & \cdots \ot x_{-1} \ot x_0 \ot \cdots.   \\
\end{array}
$$
Now, we define the linear map
$$
\begin{array}{cccc}
 U:   &   M   &  \longrightarrow   &   M  \\
      &    y &  \longmapsto       &  d\big((Id_{B(\ell^2_I)} \ot \mathscr{S})(y)\big)d.  \\
\end{array}
$$
The map $M\xra{U}M$ is a unital $*$-automorphism of $M$. Moreover,
it is easy to see that $M\xra{U}M$ preserves the trace $\tau_M$.
Now, we will show that, for any positive integer $k$, we have, for
all $x\in B\big(\ell^2_I\big)$
\begin{align}
% \nonumber to remove numbering (before each equation)
\label{hyp récurrence}
  U^k\circ J(x)
  &=\sum_{i,j\in I}x_{ij}e_{ij} \ot \cdots \ot I \ot \underbrace{\omega(e_i) \omega(e_j) \ot \cdots \ot \omega(e_i) \omega(e_j)}_{k\text{ factors}} \ot I \ot \cdots
\end{align}
by induction on $k$, where the first $\omega(e_i) \omega(e_j)$ is in
position $0$. The statement clearly holds for $k=0$. Now, assume
(\ref{hyp récurrence}). For all $x\in B\big(\ell^2_I\big)$, we have
\begin{eqnarray*}
% \nonumber to remove numbering (before each equation)
\lefteqn{U^{k+1}\circ J(x)= d\Big((Id_{B(\ell^2_I)} \ot \mathscr{S})\big(U^{k}\circ J(x)\big)\Big)d}  \\
   &=& d\Bigg((Id_{B(\ell^2_I)} \ot \mathscr{S})\Bigg(\sum_{i,j\in I}x_{ij}e_{ij} \ot \cdots \ot I\ot \omega(e_i) \omega(e_j)\ot\cdots \ot \omega(e_i) \omega(e_j)\ot I \ot \cdots\Bigg)\Bigg)d  \\
   &=& \Bigg(\sum_{r\in I} e_{rr} \ot \cdots \ot I \ot \omega(e_r) \ot I \ot \cdots\Bigg)\Bigg( \sum_{i,j\in I}x_{ij}e_{ij} \ot \cdots \ot I \ot I \ot\omega(e_i) \omega(e_j)\ot\cdots\\
   &&  \ot \omega(e_i) \omega(e_j)\ot I \ot \cdots\Bigg)\Bigg( \sum_{s\in I} e_{ss} \ot \cdots \ot I \ot \omega(e_s) \ot I \ot \cdots\Bigg) \\
   &=& \sum_{i,j,r,s\in I} x_{ij}e_{rr}e_{ij}e_{ss}  \ot \cdots \ot I \ot \omega(e_r)\omega(e_s)\ot \omega(e_i)\omega(e_j)\ot\cdots \ot \omega(e_i) \omega(e_j)\ot I \ot \cdots\\
   &=& \sum_{i,j\in I} x_{ij}e_{ij} \ot \cdots \ot I \ot \omega(e_i)\omega(e_j)\ot \omega(e_i) \omega(e_j)\ot \cdots \ot \omega(e_i) \omega(e_j)\ot I \ot \cdots.
\end{eqnarray*}
We obtained the statement (\ref{hyp récurrence}) for $k+1$. Then, we
deduce that for any positive integer $k$ and any $x\in
B\big(\ell^2_I\big)$ we have
%\pagebreak[4]
\begin{eqnarray*}
% \nonumber to remove numbering (before each equation)
 \lefteqn{ \mathbb{E}\circ U^k\circ J(x)}\\
 &=& (Id_{S^p_{I}} \ot \cdots \ot \tau \ot\cdots)\Bigg(\sum_{i,j\in I}x_{ij}e_{ij} \ot \cdots \ot I\ot \omega(e_i) \omega(e_j)\ot\cdots\ot \omega(e_i) \omega(e_j)\ot I \ot \cdots\Bigg)\\
   &=& \sum_{i,j\in I}\tau\big(\omega(e_i) \omega(e_j)\big)^k x_{ij}e_{ij}  \\
   &=& \sum_{i,j\in I}\big(\langle e_i,e_j\rangle_{\ell^{2,A}}\big)^k x_{ij}e_{ij}  \hspace{1cm} \text{by (\ref{petite formule de Wick})}\\
   &=& \sum_{i,j\in I} (a_{ij})^k x_{ij}e_{ij} \hspace{1cm} \text{by (\ref{eiej=aij})}\\
   &= &(M_A)^k(x).
\end{eqnarray*}
Thus, for any positive integer $k$, we have
\begin{equation*}
(M_A)^k=\mathbb{E}\circ U^k \circ J.
\end{equation*}
The proof is complete.
\end{preuve}

%Now we use a 2-2 matrix argument in order to show the following
%corollary.
\begin{cor} \label{coro argument matrix22} Let $B\big(\ell^2_I\big)
\xra{M_A} B\big(\ell^2_I\big)$ be a contractive Schur multiplier
associated with a real-valued matrix $A$. Suppose $1 < p < \infty$.
Then, the induced Schur multiplier $S^p_I \xra{M_A} S^p_I$ satisfies
the noncommutative Matsaev's property. More precisely, for any
complex polynomial $P$, we have
\begin{equation*}
\bnorm{P(M_A)}_{cb,S^p_I \xra{} S^p_I} \leq \norm{P}_{p,S^p}.
\end{equation*}
\end{cor}

\begin{preuve}
Suppose that $B\big(\ell^2_I\big) \xra{M_A} B\big(\ell^2_I\big)$ is
a contractive Schur multiplier associated with a real matrix $A$ of
$\mathbb{M}_I$. There exists a set $K$ and norm 1 vectors $h_i \in
\ell^2_K(\R)$ and $k_i \in \ell^2_K(\R)$ such that for all $i,j\in
I$ we have $a_{i,j}=\langle h_i,k_j\rangle_{\ell^2_K(\R)}$. Now we
define the following matrices of $\mathbb{M}_I$
$$
B=\Big[\langle h_i,h_j\rangle_{\ell^2_K(\R)}\Big]_{i,j \in I}, \ \
C=\Big[\langle k_i,k_j\rangle_{\ell^2_K(\R)}\Big]_{i,j \in I}  \ \
\text{and} \ \ D=\Big[\langle
k_i,h_j\rangle_{\ell^2_K(\R)}\Big]_{i,j \in I}.
$$
For all $i\in I$ and all $n\in \{1,2\}$, we define the norm 1 vector
$l_{(n,i)}$ of $\ell^2_{K}(\R)$ by
$$
l_{(n,i)}=\left\{
\begin{array}{cl}
      h_i    & \text{if}\  n=1   \ \text{and}\  i \in I       \\
      k_i    & \text{if}\  n=2   \ \text{and}\  i \in I.        \\
\end{array}\right.
$$
Now, by the identification $\mathbb{M}_2(\mathbb{M}_I)\simeq
\mathbb{M}_{\{1,2\}\times I}$, the matrix $ \left[
  \begin{array}{cc}
    B & A \\
    D & C \\
  \end{array}
\right] $ of $\mathbb{M}_2(\mathbb{M}_I)$ identifies
%\pagebreak[4]
to the matrix
$$
F=\Big[\langle l_{n,i},l_{m,j}\rangle_{\ell^2_{
K}(\R)}\Big]_{(n,i)\in\{1,2\}\times I ,(m,j)\in\{1,2\}\times I}
$$
of $\mathbb{M}_{\{1,2\}\times I}$. The Schur multiplier
$B\big(\ell^2_{\{1,2\} \times I}\big)\xra{M_F}B\big(\ell^2_{\{1,2\}
\times I}\big)$ associated with this matrix is unital and completely
positive. Moreover, since the matrix $F$ is real,
$B\big(\ell^2_{\{1,2\} \times I}\big)\xra{M_F}B\big(\ell^2_{\{1,2\}
\times I}\big)$ is selfadjoint. Let $1 < p < \infty$. For any
complex polynomial $P$, we have
\begin{align*}
% \nonumber to remove numbering (before each equation)
 \bnorm{P(M_A)}_{S^p_I \xra{} S^p_I }
 &\leq \left|\left|\left[
  \begin{array}{cc}
    P(M_{B}) & P(M_A) \\
    P(M_D) & P(M_{C}) \\
  \end{array}
\right] \right|\right|_{S^p_2(S^p_{I}) \xra{}S^p_2(S^p_{I}) }\\
   &=  \bnorm{P(M_F)}_{S^p_{\{1,2\} \times I} \xra{}
   S^p_{\{1,2\} \times I}}.
\end{align*}
Now, according to Theorem \ref{th dilatation SpI}, remarks following
Lemma \ref{lemma trace preserving} and Corollary \ref{dilationQWEP},
the Schur multiplier $M_F$ satisfies the noncommutative Matsaev's
property. We deduce that $M_A$ also satisfies this property.
Moreover, in applying this result to the Schur multiplier $M_{I\ot
A}(=I \ot M_{A})$, we obtain the inequality for the completely
bounded norm.
\end{preuve}

We pass to Fourier multipliers on discrete groups. Suppose that $G$
is a discrete group. We denote by $e_G$ the neutral element of $G$.
We denote by $\ell^2_G\xra{\lambda(g)}\ell^2_G$ the unitary operator
of left translation by $g$ and $\VN(G)$ the von Neumann algebra of
$G$ spanned by the $\lambda(g)$ where $g\in G$. It is an finite
algebra with trace given by
$$
\tau_{G}(x)=\big\langle\varepsilon_{e_G},x(\varepsilon_{e_G})\big\rangle_{\ell^2_G}
$$
where $(\varepsilon_g)_{g\in G}$ is the canonical basis of
$\ell^2_G$ and $x \in \VN(G)$. A Fourier multiplier is a normal
linear map $\VN(G) \xra{T} \VN(G)$ such that there exists a function
$G\xra{t}\C$ such that for all $g\in G$ we have
$T\big((\lambda(g)\big)=t_g\lambda(g)$. In this case, we denote $T$
by
$$
\begin{array}{cccc}
   M_t:   &    \VN(G)      &  \longrightarrow   & \VN(G)   \\
          &    \lambda(g)  &  \longmapsto       & t_g\lambda(g).   \\
\end{array}
$$
It is easy to see that a contractive Fourier multiplier $\VN(G)
\xra{M_t} \VN(G)$ is selfadjoint if and only if $G\xra{t}\C$ is a
real function. It is well-known that a Fourier multiplier $\VN(G)
\xra{M_t} \VN(G)$ is completely positive if and only if the function
$t$ is positive definite. If the discrete group $G$ is amenable, by
\cite[Corollary 1.8]{DCH}, every contractive Fourier multiplier
$\VN(G) \xra{M_t} \VN(G)$ is completely contractive. Recall the
following particular case of the transfer result \cite[Theorem
2.7]{NeR}.
\begin{thm}
Let $G$ be an amenable discrete group. Suppose $1<p<\infty$. Let $G
\xra{t}\R$ a function. Let $A$ be the matrix of $\mathbb{M}_{G}$
defined by $a_{g,h}=t_{gh^{-1}}$ where $g,h \in G$. The Fourier
multiplier $M_t$ is completely bounded on $L^p(\VN(G))$ if and only
if the Schur multiplier $M_A$ is completely bounded on $S^p_{G}$. In
this case, we have
\begin{equation}\label{transfer=}
\bnorm{M_t}_{cb,L^p(\VN(G)) \xra{}
L^p(\VN(G))}=\bnorm{M_A}_{cb,S^p_{G}\xra{}S^p_{G}}.
\end{equation}
\end{thm}
Now, we can prove the next result.
\begin{cor}
Let $G$ be an amenable discrete group. Let $\VN(G) \xra{M_t} \VN(G)$
be a contractive Fourier multiplier associated with a real function
$G \xra{t}\R$. Suppose $1 < p < \infty$. Then, the induced Fourier
multiplier $L^p\big(\VN(G)\big) \xra{M_t}L^p\big(\VN(G)\big)$
satisfies the noncommutative Matsaev's property. More precisely, for
any complex polynomial $P$, we have
$$
\bnorm{P(M_t)}_{cb,L^p(\VN(G)) \xra{} L^p(\VN(G))} \leq
\norm{P}_{p,S^p}.
$$
\end{cor}

\begin{preuve}
We define the matrix $A$ of $\mathbb{M}_{G}$ by
$a_{g,h}=t_{gh^{-1}}$ where $g,h \in G$. By (\ref{transfer=}), for
any complex polynomial $P$ and all $1 < p < \infty$, we have
\begin{eqnarray*}
% \nonumber to remove numbering (before each equation)
\bnorm{P(M_t)}_{cb,L^p(\VN(G)) \xra{}L^p(\VN(G))}
   &=& \bnorm{P(M_A)}_{cb,S^p_{G}\xra{}S^p_{G}}  \\
   &=& \bnorm{P(Id_{S^p} \ot M_A)}_{S^p(S^p_{G}) \xra{}S^p(S^p_{G})}  \\
   &=& \bnorm{P(M_{I \ot A})}_{S^p(S^p_{G}) \xra{}S^p(S^p_{G})}.
\end{eqnarray*}
Since $G$ is amenable, the Fourier multiplier $\VN(G) \xra{M_t}
\VN(G)$ is completely contractive. Moreover, the map $\VN(G)
\xra{M_t} \VN(G)$ is selfadjoint. Thus, for any $1 < p < \infty$,
the map $L^p\big(\VN(G)\big) \xra{M_t}L^p\big(\VN(G)\big)$ is
completely contractive. Then, by (\ref{transfer=}), for any $1 < p <
\infty$, we have
\begin{eqnarray*}
% \nonumber to remove numbering (before each equation)
\norm{M_{I\ot A }}_{S^p(S^p_{G}) \xra{}S^p(S^p_{G})}
   &=&  \norm{M_A}_{cb,S^p_{G}\xra{} S^p_{G}} \\
   &=& \norm{M_t}_{cb,L^p(\VN(G)) \xra{}L^p(\VN(G))} \\
   &\leq& 1.
\end{eqnarray*}
By Corollary 4.3, we deduce finally that, for any complex polynomial
$P$ and all $1 < p < \infty$, we have
$$
\bnorm{P(M_t)}_{L^p(\VN(G)) \xra{}L^p(\VN(G))} \leq
\norm{P}_{p,S^p(S^p_{G})}=\norm{P}_{p,S^p}.
$$
\end{preuve}

In order to prove the next theorem we need the notion of
crossproduct. We refer to \cite{Stra} and \cite{Sun} for more
information on this concept.
\begin{thm}
\label{th dilatation Fourier multipliers2} Let $G$ be a discrete
group. Let $\VN(G) \xra{M_t} \VN(G)$ be a unital completely positive
Fourier multiplier associated with a real valued function $G
\xra{t}\R$. Then there exists a von Neumann algebra $M$ equipped
with a faithful finite normal trace, a  unital trace preserving
$*$-automorphism $M \xra{U} M$, a unital normal trace preserving
one-to-one $*$-homomorphism $\VN(G) \xra{J} M$ such that,
\begin{eqnarray*}
(M_t)^k=\mathbb{E}U^kJ
\end{eqnarray*}
for any integer $k\geq 0$, where $M \xra{\mathbb{E}} VN(G)$ is the
canonical faithful normal trace preserving conditional expectation
associated with $J$.
\end{thm}

\begin{preuve}
Since the map $\VN(G) \xra{M_t} \VN(G)$ is completely positive, we
can define a positive symmetric bilinear form
$\langle\cdot,\cdot\rangle_{\ell^{2,t}}$ on the real span of the
$e_g$, where $g \in G$, by:
$$
\langle e_g,e_h\rangle_{\ell^{2,t}}=t_{g^{-1}h}.
$$
We denote by $\ell^{2,t}$ the completion of the real pre-Hilbert
space obtained by quotient by the corresponding kernel. For all
$g\in G$, we denote by $g$ the class of $e_g$ in $\ell^{2,t}$. Then,
for all $g,h\in G$, we have
$$
\langle g,h\rangle_{\ell^{2,t}}=t_{g^{-1}h}.
$$
For all $g\in G$, it easy to see that there exists a unique isometry
$\ell^{2,t}\xra{\theta_g}\ell^{2,t}$ such that for all $ h\in G$ we
have $\theta_g(h)=gh.$

For all $g\in G$, we define the unital trace preserving
$*$-automorphism $\alpha(g)=\Gamma_{-1}(\theta_g \ot
Id_{\ell^{2}_{\mathbb{Z}}})$:
$$
\begin{array}{cccc}
   \alpha(g): &  \Gamma_{-1}\big(\ell^{2,t}\ot_2 \ell^{2}_{\mathbb{Z}}\big)   &  \longrightarrow   &  \Gamma_{-1}\big(\ell^{2,t}\ot_2 \ell^{2}_{\mathbb{Z}}\big)  \\
              &   w(h \ot v)  &  \longmapsto       &  w(gh \ot v).  \\
\end{array}
$$
The homomorphism
$G\xra{\alpha}\Aut\Big(\Gamma_{-1}\big(\ell^{2,t}\ot_2
\ell^{2}_{\mathbb{Z}}\big)\Big)$ allows us to define the von Neumann
algebra
\begin{eqnarray}
% \nonumber to remove numbering (before each equation)
\label{von Neumann algebra M}
  M &=&\Gamma_{-1}\big(\ell^{2,t}\ot_2\ell^{2}_{\mathbb{Z}}\big)\rtimes_\alpha G.
\end{eqnarray}
We can identify
$\Gamma_{-1}\big(\ell^{2,t}\ot_2\ell^{2}_{\mathbb{Z}}\big)$ as a
subalgebra of $M$. We let $J$ the canonical normal unital injective
$*$-homomorphism $\VN(G) \xra{J} M$. We denote by $\tau$ the
faithful finite normal trace on
$\Gamma_{-1}\big(\ell^{2,t}\ot_2\ell^{2}_{\mathbb{Z}}\big)$. Recall
that, for all $g\in G$, the map $\alpha(g)$ is trace preserving.
Thus, the trace $\tau$ is $\alpha$-invariant. We equip $M$ with the
induced canonical trace $\tau_M$. For all $x\in
\Gamma_{-1}\big(\ell^{2,t}\ot_2\ell^{2}_{\mathbb{Z}}\big)$ and all
$g\in G$, we have
\begin{equation}\label{formule trace}
\tau_M\Big(xJ\big(\lambda(g)\big)\Big)=\delta_{g,e_G}\tau(x)
\end{equation}
(see \cite{Stra} pages 359 and 352). If $g,h\in G$ and $v\in
\ell^{2}_{\mathbb{Z}}$, we can write the relations of commutation of
the crossed product as
%\pagebreak[4]
\begin{equation} \label{commutation crossed
product} J\big(\lambda(g)\big)\omega(h \ot v)=\omega(gh \ot
v)J\big(\lambda(g)\big).
\end{equation}
We denote by $M \xra{\mathbb{E}} \VN(G)$ the canonical faithful
normal trace preserving conditional expectation. For all $x\in
\Gamma_{-1}\big(\ell^{2,t}\ot_2 \ell^{2}_{\mathbb{Z}}\big)$ and all
$g\in G$ we have
$$
\mathbb{E}\Big(xJ\big(\lambda(g)\big)\Big)=\tau(x)\lambda(g).
$$
We define the unital trace preserving $*$-automorphism
$\mathscr{S}=\Gamma_{-1}(Id_{\ell^{2,t}}\ot S)$:
$$
\begin{array}{cccc}
  \mathscr{S}:      &  \Gamma_{-1}\big(\ell^{2,t}\ot_2 \ell^{2}_{\mathbb{Z}}\big)   &  \longrightarrow   &  \Gamma_{-1}\big(\ell^{2,t}\ot_2 \ell^{2}_{\mathbb{Z}}\big)  \\
                    &                \omega(h\ot e_n)                               &  \longmapsto       &  \omega(h\ot e_{n+1}).   \\
\end{array}
$$
Since $M_t$ is unital, $\omega(e_G\ot e_0)$ is a symmetry, i.e. a
selfadjoint unitary element. Moreover, for all $g\in G$, we have
$\alpha(g)\mathscr{S}=\mathscr{S}\alpha(g)$. Then, by
\cite[Proposition 4.4.4]{Sun}, we can define a unital
$*$-automorphism
$$
\begin{array}{cccc}
  U:      &                  M             &  \longrightarrow   &  M  \\
          &   xJ\big(\lambda(g)\big)       &  \longmapsto       & \omega(e_G \ot e_0)\mathscr{S}(x)J\big(\lambda(g)\big)\omega(e_G \ot e_0).   \\
\end{array}
$$
Now, we will show that $U$ preserves the trace. For all $g \in G$
and all $x \in
\Gamma_{-1}\big(\ell^{2,t}\ot_2\ell^{2}_{\mathbb{Z}}\big)$, we have
\begin{eqnarray*}
% \nonumber to remove numbering (before each equation)
\tau_M\bigg(U\Big(xJ\big(\lambda(g)\big)\Big)\bigg)
   &=& \tau_M\Big(\mathscr{S}(x)J\big(\lambda(g)\big)\Big)  \\
   &=& \delta_{g,e_G}\tau\big(\mathscr{S}(x)\big)\hspace{1cm} \text{by (\ref{formule trace})}\\
   &=& \delta_{g,e_G}\tau(x)\\
   &=& \tau_M\Big(xJ\big(\lambda(g)\big)\Big).
\end{eqnarray*}
We conclude by linearity and normality. It is not hard to see that
$J$ preserves the traces.

Now, we will prove that, for any integer $k\geq 1$ and any $g\in G$,
we have
\begin{eqnarray}
\label{hypo recurrence 2}
% \nonumber to remove numbering (before each equation)
 \lefteqn{U^k\circ J\big(\lambda(g)\big)} && \\
 &=& \omega(e_G \ot e_0)\omega(e_G \ot e_1)\cdots \omega(e_G \ot e_{k-1}) \omega(g \ot e_{k-1})\omega(g \ot e_{k-2})\cdots \omega(g \ot e_0)J\big(\lambda(g)\big) \nonumber
\end{eqnarray}
by induction on $k$. The statement holds clearly for $k=1$: if $g\in
G$, we have %\pagebreak[4]
\begin{eqnarray*}
% \nonumber to remove numbering (before each equation)
 U\circ J\big(\lambda(g)\big)
   &=&  \omega(e_G \ot e_0)J\big(\lambda(g)\big)\omega(e_G \ot e_0)\\
   &=& \omega(e_G \ot e_0)\omega(g \ot e_0)J\big(\lambda(g)\big) \hspace{1cm} \text{by (\ref{commutation crossed product})}.
\end{eqnarray*}
Now, assume (\ref{hypo recurrence 2}). For all $g\in G$, we have
\begin{eqnarray*}
% \nonumber to remove numbering (before each equation)
 \lefteqn{ U^{k+1}\circ J\big(\lambda(g)\big) }\\
 &=& U\Big(\omega(e_G \ot e_0)\omega(e_G \ot e_1)\cdots \omega(e_G \ot e_{k-1}) \omega(g \ot e_{k-1})\omega(g \ot e_{k-2})\cdots  \omega(g \ot e_0)J\big(\lambda(g)\big)\Big)  \\
 &=& \omega(e_G \ot e_0)\omega(e_G \ot e_1)\omega(e_G \ot e_2)\cdots \omega(e_G \ot e_{k}) \omega(g \ot e_{k})\omega(g \ot e_{k-1})\cdots   \\
 &&  \omega(g \ot e_1)J\big(\lambda(g)\big)\omega(e_G \ot e_0)\\
 &=& \omega(e_G \ot e_0)\omega(e_G \ot e_1)\omega(e_G \ot e_2)\cdots \omega(e_G \ot e_{k}) \omega(g \ot e_{k})\omega(g \ot e_{k-1})\cdots  \\
 && \omega(g \ot e_1)\omega(g \ot e_0)J\big(\lambda(g)\big).
\end{eqnarray*}
We obtained the statement (\ref{hypo recurrence 2}) for $k+1$. Now
let $k\geq 1$ and $g\in G$. We define the elements $f_1,\ldots
,f_{2k}$ of the Hilbert space $\ell^{2,t}
\ot_2\ell^{2}_{\mathbb{Z}}$ by
$$
f_i = \left\{
\begin{array}{cl}
   e_G \ot e_{i-1} & \text{if}\ 1 \leq i\leq k    \\
   g \ot e_{2k-i}    &  \text{if}\  k+1 \leq i\leq 2k  \\
\end{array}
\right.
$$
If $ 1\leq i\leq 2k$, we have
\begin{eqnarray*}
% \nonumber to remove numbering (before each equation)
\langle f_{i},f_{2k-i+1}\rangle_{\ell^{2,t}
\ot_2\ell^{2}_{\mathbb{Z}}} &=&
\langle e_G \ot e_{i-1},g \ot e_{i-1}\rangle_{\ell^{2,t}\ot_2\ell^{2}_{\mathbb{Z}}} \\
   &=& \langle e_G,g \rangle_{\ell^{2,t}}\langle e_{i-1},e_{i-1}\rangle_{\ell^{2}_{\mathbb{Z}}} \\
   &=& t_g
\end{eqnarray*}
By a similar computation, if $ 1\leq i < j\leq 2k$ with
$j\not=2k-i+1$, we obtain $\langle f_i,f_j\rangle_{\ell^{2,t}
\ot_2\ell^{2}_{\mathbb{Z}}}=0$. Then, for all $g\in G$, we have
\begin{eqnarray}
% \nonumber to remove numbering (before each equation)
\mathbb{E} U^k J\big(\lambda(g)\big) %\nonumber\\ \lefteqn{
 &=& \mathbb{E}\Big( \omega(e_G \ot e_0)\cdots \omega(e_G \ot e_{k-1}) \omega(g \ot e_{k-1})\cdots  \omega(g \ot e_0)J\big(\lambda(g)\big)\Big) \nonumber \\
 &=& \tau\big( \omega(e_G \ot e_0)\cdots \omega(e_G \ot e_{k-1}) \omega(g \ot e_{k-1})\cdots  \omega(g \ot e_0)\big)\lambda(g)  \nonumber\\
 &=& \tau\big( \omega(f_1)\omega(f_2)\cdots \omega(f_{2k})\big)\lambda(g)  \nonumber\\
 &=& \Bigg(\sum_{\mathcal{V}\in \mathcal{P}_2(2k)} (-1)^{c(\mathcal{V})} \prod_{(i,j)\in \mathcal{V}}\langle f_i,f_j\rangle_{\ell^{2,t}\ot_2\ell^{2}_{\mathbb{Z}}}\Bigg)\lambda(g)\hspace{1cm} \text{by (\ref{formule de Wick})}\label{Fock product in dilatation}\\
 &=& \langle f_1,f_{2k}\rangle_{\ell^{2,t} \ot_2\ell^{2}_{\mathbb{Z}}}\cdots \langle f_k,f_{k+1}\rangle_{\ell^{2,t} \ot_2\ell^{2}_{\mathbb{Z}}}\lambda(g)\nonumber\\
%% &=& <e_G \ot e_0,g \ot e_0>_{\ell^{2,t} \ot_2\ell^{2}_{\mathbb{Z}}}...<e_G \ot e_{k-1},g \ot e_{k-1}>_{\ell^{2,t} \ot_2\ell^{2}_{\mathbb{Z}}}\lambda(g)\\
%% &=& <e_G,g>_{\ell^{2,t}}<e_0,e_0>_{\ell^{2}_{\mathbb{Z}}}...<e_G,g>_{\ell^{2,t}}<e_{k-1},e_{k-1}>_{\ell^{2}_{\mathbb{Z}}}\lambda(g)\nonumber\\
 &=& (t_{g})^k\lambda(g)\nonumber\\
 &=& T^k\big(\lambda(g)\big).\nonumber
\end{eqnarray}
$\big($The only non-zero term in the sum of (\ref{Fock product in
dilatation}) is the term with $
\mathcal{V}=\big\{(1,2k),(2,2k-1),\ldots,(k,k+1)\big\}$, which
satisfies $c(\mathcal{V})=0\big)$. Thus, for any positive integer
$k$ (the case $k=0$ is trivial), we conclude that
$$
T^k=\mathbb{E} U^k J.
$$
\end{preuve}

\begin{cor}
Let $G$ be a discrete group. Let $\VN(G) \xra{M_t} \VN(G)$ be a
unital completely positive Fourier multiplier which is associated
with a real function $G \xra{t}\R$. Suppose that the von Neumann
algebra (\ref{von Neumann algebra M}) has QWEP. Let $1 \leq p \leq
\infty$. Then, the induced Fourier multiplier $L^p\big(\VN(G)\big)
\xra{M_t}L^p\big(\VN(G)\big)$ satisfies the noncommutative Matsaev's
property.
\end{cor}

\begin{preuve}
This corollary follows from Theorem \ref{th dilatation Fourier
multipliers2}, remarks following Lemma \ref{lemma trace preserving}
and Corollary \ref{dilationQWEP}.
\end{preuve}

At the light of above corollary, it is important to know when the
von Neumann algebra (\ref{von Neumann algebra M}) has QWEP. If the
group $G$ is amenable, this algebra has QWEP by \cite[Proposition
4.1]{Oza} (or \cite[Proposition 6.8]{Con}). Now we give an example
of non-amenable group $G$ such that this von Neumann algebra has
QWEP. We denote by $\mathbb{F}_n$ a free group with $n$ generators
denoted by $g_1,\ldots,g_n$ where $1\leq n \leq \infty$. We denote
by $\mathcal{R}$ the hyperfinite factor of type ${\rm II_1}$ and by
$\mathcal{R}^\ull$ an ultrapower of $\mathcal{R}$ with respect to a
non-trivial ultrafilter $\ull$. In order to prove the next theorem
we need the notion of amalgamated free product of von Neumann
algebras. We refer to \cite{BlD} and \cite{Ued} for more information
on this concept. Note that, with the notations of the proof of
Theorem \ref{th dilatation Fourier multipliers2}, the von Neumann
algebra $\Gamma_{-1}\big(\ell^{2,t}\ot_2\ell^{2}_{\Z}\big)$ is
$*$-isomorphic to the hyperfinite factor of type $\rm{II}_1$.
\begin{prop}
\label{prop AVN has QWEP} Suppose $1 \leq n \leq \infty$. Let
$\mathbb{F}_n \xra{\alpha} \Aut(\mathcal{R})$ be a homomorphism.
Then the crossed product $\mathcal{R}\rtimes_\alpha \mathbb{F}_n$
has QWEP.
\end{prop}

\begin{preuve}
First we will show the result for $n=2$. We denote by $\langle
g_1\rangle$ and $\langle g_2\rangle$ the subgroups of $\mathbb{F}_2$
generated by $g_1$ and $g_2$ and by $\alpha_1$ and $\alpha_2$ the
restrictions of $\alpha$ to these subgroups. First, we prove that
the subalgebras $\mathcal{R} \rtimes_{\alpha_1} \langle g_1\rangle$
and $\mathcal{R}\rtimes_{\alpha_2} \langle g_2\rangle$ are free with
respect to the canonical faithful normal trace preserving
conditional expectation $\mathcal{R}\rtimes_\alpha
\mathbb{F}_2\xra{\mathbb{E}}\mathcal{R}$. We identify $\mathcal{R}$
as a subalgebra of $\mathcal{R}\rtimes_\alpha \mathbb{F}_2$. We may
regard the elements of $\mathcal{R}\rtimes_\alpha \mathbb{F}_2$ as
matrices $\Big[\alpha_{r^{-1}}\big(\varpi(rt^{-1})\big)\Big]_{r,t\in
\mathbb{F}_2}$ with entries in $\mathcal{R}$ where
$\mathbb{F}_2\xra{\varpi}\mathcal{R}$ is a map. Recall that the
conditional expectation $\mathbb{E}$ on $\mathcal{R}$ is given by
$$
\mathbb{E}\bigg(\Big[\alpha_{r^{-1}}\big(\varpi(rt^{-1})\big)\Big]_{r,t\in
\mathbb{F}_2}\bigg)=\varpi(e_{\mathbb{F}_2}).
$$
Suppose that $i_1,\ldots,i_k\in \{1,2\}$ are integers such that
$i_1\not=i_2,\ldots,i_{k-1}\not=i_{k}$. For any $1 \leq j \leq k$,
let
$$
A_{j}=\Big[\alpha_{r^{-1}}\big(\varpi_{j}(rt^{-1})\big)\Big]_{r,t\in
\mathbb{F}_2}
$$
be an element of $\mathcal{R}\rtimes_{\alpha_{i_j}} \langle
g_{i_j}\rangle $ such that $\mathbb{E}(A_{j})=0$ where each
$\mathbb{F}_2\xra{\varpi_{j}}\mathcal{R}$ is a map satisfying
$\varpi_{j}(g)=0$ if $g \not\in \langle g_{i_j}\rangle$. Then, for
all $1 \leq j \leq k$, we have $\varpi_{j}(e_{\mathbb{F}_2})=0$.
Now, we have
\begin{eqnarray*}
% \nonumber to remove numbering (before each equation)
\mathbb{E}(A_1\cdots A_k)
   &=& \mathbb{E}\bigg(\Big[\alpha_{r^{-1}}\big(\varpi_{1}(rt^{-1})\big)\Big]_{r,t\in \mathbb{F}_2}\cdots\Big[\alpha_{r^{-1}}\big(\varpi_{k}(rt^{-1})\big)\Big]_{r,t\in \mathbb{F}_2}\bigg)\\
   &=& \sum_{L^1,...,l_{k-1} \in
   \mathbb{F}_2}\varpi_{1}\big(l_{1}^{-1}\big)\alpha_{l_1^{-1}}\Big(\varpi_{2}\big(l_{1}l_{2}^{-1}\big)\Big)\cdots\alpha_{l_{k-2}^{-1}}\Big(\varpi_{2}\big(l_{k-2}l_{k-1}^{-1}\big)\Big)\alpha_{l_{k-1}^{-1}}\big(\varpi_{k}(l_{k-1})\big)\\
   &=& 0.
\end{eqnarray*}
%Note that, by hypothesis, we have $\varpi_{1}\big(l_{1}^{-1}\big)=0$
%if $L^1^{-1} \not\in \langle g_{i_1}\rangle$,
%$\varpi_{2}\big(l_{1}l_{2}^{-1}\big)=0$ if $l_{1}l_{2}^{-1}\in
%\langle g_{i_2}\rangle$,..., $\varpi_{k}(l_{k-1})=0$ if $l_{k-1}\in
%\langle g_{i_k}\rangle$. Then, it is not difficult to see that
%$$
%\mathbb{E}(A_1\cdots A_k)=0.
%$$                                                              embeddable into $\mathcal{R}^\ull$.
Thus the von Neumann algebra $\mathcal{R}\rtimes_\alpha
\mathbb{F}_2$ decomposes as an amalgamated free product of
$\mathcal{R}\rtimes_{\alpha_1} \langle g_1 \rangle$ and
$\mathcal{R}\rtimes_{\alpha_2} \langle g_2\rangle$ over
$\mathcal{R}$. Moreover, the groups $\langle g_1 \rangle$ and
$\langle g_2 \rangle$ are commutative, hence amenable. We have
already point out that the crossed product of the hyperfinite factor
$\mathcal{R}$ by an amenable group has QWEP. Then the von Neumann
algebras $\mathcal{R}\rtimes_{\alpha_1} \langle g_1 \rangle$ and
$\mathcal{R}\rtimes_{\alpha_2} \langle g_2\rangle$ are QWEP.
Moreover, by \cite[page 283]{Bla}, these von Neumann algebras have a
separable predual. By \cite[Theorem 1.4]{Kir}, we deduce that these
von Neumann algebras are embeddable into $\mathcal{R}^\ull$. Now,
the theorem stating in \cite[Corollary 4.5]{BDJ} says that, for
finite von Neumann algebras with separable preduals, being
embeddable into $\mathcal{R}^\ull$ is stable under amalgamated free
products over a hyperfinite von Neumann algebra. Thus we deduce that
$\mathcal{R}\rtimes_\alpha \mathbb{F}_2$ is embeddable into
$\mathcal{R}^\ull$, which is equivalent to QWEP, by \cite[Theorem
1.4]{Kir}, since $\mathcal{R}\rtimes_\alpha \mathbb{F}_2$ has a
separable predual. Induction then gives the case when $2 \leq n
<\infty$, and the case $n=\infty$ then follows since, by
\cite[Proposition 4.1]{Oza}, QWEP is preserved by taking the weak*
closure of increasing unions of von Neumann algebras.
\end{preuve}

We pass to maps arising in the second quantization in the context of
\cite{BKS}.
\begin{prop}
Suppose $1 < p < \infty$ and $-1\leq q < 1$. Let $H$ be a real
Hilbert space and $H \xra{T} H$ a contraction. Then the induced map
$L^p\big(\Gamma_q(H)\big) \xra{M_t}L^p\big(\Gamma_q(H)\big)$
satisfies the noncommutative Matsaev's property.
\end{prop}

\begin{preuve}
There exists an orthogonal dilation $K\xra{U}K$ de $H \xra{T} H$. We
denote by $H \xra{J} K$ the embedding of $H$ in $K$ and $K \xra{Q}
H$ the projection of $K$ on $H$. The map $\Gamma_q(K)
\xra{\Gamma(J)} \Gamma_q(K)$ is a unital injective normal trace
preserving $*$-homomorphism. The map $\Gamma_q(H) \xra{\Gamma(U)}
\Gamma_q(K)$ is a unital trace preserving $*$-automorphism. The map
$\Gamma_q(K) \xra{\Gamma(Q)} \Gamma_q(H)$ is the canonical faithful
normal unital trace preserving conditional expectation of
$\Gamma_q(K)$ on $\Gamma_q(H)$. Moreover, we have for any integer
$k$
$$
\Gamma_q(T)^k =\Gamma_q(P)\Gamma_q(U)^k\Gamma_q(J).
$$
We conclude with Theorem \ref{th dilatation Fourier multipliers2},
remarks following Lemma \ref{lemma trace preserving}, Corollary
\ref{dilationQWEP} and by using the fact that, by \cite{Nou}, the
von Neumann algebra $\Gamma_q(H)$ has QWEP.
\end{preuve}

In order to state more easily our following result we need to define
the following property.  Let $M$ be a von Neumann algebra. Suppose
that $M\xra{T}M$ is a linear map.

\begin{Property}
\label{property} There exists a von Neumann algebra $N$ with QWEP
equipped with a normal faithful finite trace on $N$, a unital trace
preserving $*$-automorphism $N \xra{U} N$, a unital injective normal
trace preserving $*$-homomorphism $M \xra{J} N$ such that,
$$
T^k=\mathbb{E}U^kJ.
$$
for any integer $k\geq 0$, where $M \xra{\mathbb{E}} VN(G)$ is the
canonical faithful normal trace preserving conditional expectation
associated with $J$.
\end{Property}

This property is stable under free product. Indeed, one can prove
the next proposition with an argument similar to that used in the
proof of \cite[Lemma 10.4]{JMX} and by using \cite[Corollary
4.5]{BDJ} and \cite[Theorem 1.4]{Kir}.

\begin{prop}
Let $M_1$ and $M_2$ be von Neumann algebras with separable preduals
equipped with normal faithful finite traces $\tau_1$ and $\tau_2$.
Let $M_1\xra{T_1}M_1$ and $M_2\xra{T_2}M_2$ be linear maps. If $T_1$
and $T_2$ satisfy Property \ref{property}, their free product
$$
(M_1,\tau_1)\overline{*}(M_2,\tau_2)\xra{T_1\overline{*}T_2}(M_1,\tau_1)\overline{*}(M_2,\tau_2)
$$
also satisfies Property \ref{property}.
\end{prop}
Thus the above proposition allows us to construct other examples of
contractions satisfying the noncommutative Matsaev's property.

%%%%%%%%%%%%%%%%%%%%%%%%%%%%%%%%%%%%%%%%%%%%%%%%%%%%%%%%%%%%%%%%
%%%%%%%%%%%%%%%%%%%%%%%%%%%%%%%%%%%%%%%%%%%%%%%%%%%%%%%%%%%%%%%%

\section{The case of semigroups}

%%%%%%%%%%%%%%%%%%%%%%%%%%%%%%%%%%%%%%%%%%%%%%%%%%%%%%%%%%%%%%%%
%%%%%%%%%%%%%%%%%%%%%%%%%%%%%%%%%%%%%%%%%%%%%%%%%%%%%%%%%%%%%%%%

Suppose $1 \leq p < \infty$. We denote by $(\mathcal{T}_t)_{t \geq
0}$ the translation semigroup on $L^p(\mathbb{R})$, where
$\mathcal{T}_t(f)(s)=f(s-t)$ if $f \in L^p(\mathbb{R})$ and $s,t\in
\R$. This semigroup $(\mathcal{T}_t)_{t \geq 0}$ is a
$C_0$-semigroup of contractions.

Let $(T_t)_{t\geq 0}$ be a $C_0$-semigroup of contractions on a
Banach space $X$. For all $b \in L^1(\mathbb{R})$ with support in
${\mathbb{R}^+}$, it is easy to see that the linear operator
$$
\begin{array}{cccc}
 \int_{0}^{+\infty}b(t)T_t dt:  &       X        &  \longrightarrow   &           X                    \\
                                &       x        &  \longmapsto       &  \int_{0}^{+\infty}b(t)T_tx dt  \\
\end{array}
$$
is well-defined and bounded. Moreover, we have
\begin{eqnarray*}
% \nonumber to remove numbering (before each equation)
  \Bgnorm{\int_{0}^{+\infty}b(t)T_t dt}_{X \xra{} X}
&\leq& \norm{b}_{L^1(\R)}.
\end{eqnarray*}

Now, let us state a question for semigroups which is analogue to
Matsaev's Conjecture \ref{Matsaev}.
\begin{quest}
\label{question semigroup} Suppose $1 < p < \infty$, $p\not=2$. Let
$(T_t)_{t\geq 0}$ be a $C_0$-semigroup of contractions on a
$L^p$-space $L^p(\Omega)$ of a measure space $\Omega$. Do we have
the following estimate
\begin{equation}
\label{inequality semigroup} \Bgnorm{\int_{0}^{+\infty}b(t)T_t
dt}_{L^p(\Omega) \xra{} L^p(\Omega)} \leq
\Bgnorm{\int_{0}^{+\infty}b(t)\mathcal{T}_t dt}_{L^p(\mathbb{R})
\xra{}L^p(\mathbb{R})}
\end{equation}
for all $b \in L^1(\mathbb{R})$ with support in ${\mathbb{R}^+}$?
\end{quest}

We pass to the noncommutative case. We can state the following
noncommutative analogue of Question \ref{question semigroup}.
\begin{quest}
\label{question NC semigroup} Suppose $1< p < \infty$, $p\not=2$.
Let $(T_t)_{t\geq 0}$ be a $C_0$-semigroup of contractions on a
noncommutative $L^p$-space $L^p(M)$. Do we have the following
estimate
\begin{equation}
\label{inequality NC semigroup} \Bgnorm{\int_{0}^{+\infty}b(t)T_t
dt}_{L^p(M) \xra{} L^p(M)} \leq
\Bgnorm{\int_{0}^{+\infty}b(t)\mathcal{T}_t dt}_{cb,L^p(\mathbb{R})
\xra{}L^p(\mathbb{R})}
\end{equation}
for all $b \in L^1(\mathbb{R})$ with support in ${\mathbb{R}^+}$?
\end{quest}

For all $b \in L^1(\mathbb{R})$ with support in ${\mathbb{R}^+}$, it
is clear that $C_b=\int_{0}^{+\infty}b(t)\mathcal{T}_t dt$.
Moreover, for all $b \in L^1(\mathbb{R})$, we have
$$
\norm{C_b}_{L^1(\R)\xra{}L^1(\R)}=\norm{C_b}_{cb,L^1(\R)\xra{}L^1(\R)}
=\norm{b}_{L^1(\R)}.
$$
Consequently, the inequalities (\ref{inequality semigroup}) and
(\ref{inequality NC semigroup}) hold true for $p=1$.

In \cite[page 25]{CoW}, it is proved that the $C_0$-semigroups of
positive contractions satisfy inequality (\ref{inequality
semigroup}). Using \cite[Theorem 3]{Pel1} and the same method, we
can generalize this result to $C_0$-semigroups of operators which
admit a contractive majorant. Now, we adapt this method in order to
give a link between Question \ref{question NC semigroup} and
Question \ref{quest matsaev semifinite}.
\begin{thm}
Suppose $1< p < \infty$. Let $(T_t)_{t\geq 0}$ be a $C_0$-semigroup
of contractions on a noncommutative $L^p$-space $L^p(M)$ such that
each $L^p(M)\xra{T_t}L^p(M)$ satisfies the noncommutative Matsaev's
property. Then the semigroup $(T_t)_{t\geq 0}$ satisfies inequality
(\ref{inequality NC semigroup}).
\end{thm}

\begin{preuve}
It is not hard to see that it suffices to prove this in the case
when $b$ has compact support. Now we define the sequence
$\big(a_n\big)_{n\geq 1}$ of complex sequences indexed by
$\mathbb{Z}$ as in the proof of Theorem \ref{mult bounded not
completely bounde}. Let $n\geq 1$. Observe that if $\mathbb{R}^+
\xra{f} L^p(M)$ is continuous and piecewise affine with nodes at
$\frac{k}{n}$ then
$$
\int_{0}^{+\infty}b(t)f(t)dt=\sum_{k=0}^{+\infty}
a_{n,k}f\bigg(\frac{k}{n}\bigg).
$$
Let $x \in L^p(M)$. Let $\mathbb{R}^+ \xra{f_n} L^p(M)$ be the
continuous and piecewise affine function with nodes at $\frac{k}{n}$
such that $ f_n\big(\frac{k}{n}\big)=\big(T_{\frac{1}{n}}\big)^k x$.
Since the map $t \mapsto T_tx$ is uniformly continuous on compacts
of $\mathbb{R}^+$ we have
\begin{eqnarray*}
% \nonumber to remove numbering (before each equation)
 \Bgnorm{\int_{0}^{+\infty}b(t)T_t xdt-\sum_{k=0}^{+\infty}a_{n,k}\big(T_{\frac{1}{n}}\big)^k x}_{L^p(M)}
   &=& \Bgnorm{\int_{0}^{+\infty}b(t)T_t xdt-\sum_{k=0}^{+\infty}a_{n,k}f_{n}\bigg(\frac{k}{n}\bigg)}_{L^p(M)} \\
   &=& \Bgnorm{\int_{0}^{+\infty} b(t)\big(T_tx-f_n(t)\big)dt}_{L^p(M)}\xra[n \to
   +\infty]{}0.
\end{eqnarray*}
We deduce that
$$
\sum_{k=0}^{+\infty} a_{n,k}\big(T_{\frac{1}{n}}\big)^{k} \xra[n \to
+\infty]{so} \int_{0}^{+\infty} b(t)T_t dt.
$$
By the commutative diagram of the proof of Theorem \ref{mult bounded
not completely bounde}, we have for any integer $n\geq 1$
$$
\bnorm{C_{a_n}}_{cb,\ell^p_\mathbb{Z} \xra{} \ell^p_\mathbb{Z}} \leq
\norm{C_b}_{cb,L^p(\mathbb{R}) \xra{} L^p(\mathbb{R})}.
$$
Finally, by the strongly lower semicontinuity of the norm, we obtain
that %\pagebreak[4]
\begin{eqnarray*}
% \nonumber to remove numbering (before each equation)
 \Bgnorm{\int_{0}^{+\infty}b(t)T_t \ dt}_{L^p(M)\xra{}L^p(M)}
   &\leq& \liminf_{n \to +\infty} \Bgnorm{\sum_{k=1}^{+\infty} a_{n,k} \big(T_{\frac{1}{n}}\big)^k}_{L^p(M)\xra{}L^p(M)} \\
   &\leq&    \liminf_{n \to +\infty} \norm{C_{a_n}}_{cb,\ell^p_\mathbb{Z} \xra{} \ell^p_\mathbb{Z}}  \\
   &=&    \norm{C_b}_{cb,L^p(\mathbb{R}) \xra{} L^p(\mathbb{R})}.
\end{eqnarray*}
\end{preuve}

The first consequence of this theorem is that inequality
(\ref{inequality NC semigroup}) holds true for $p=2$. Now, we list
some natural examples of semigroups which satisfy inequality
(\ref{inequality NC semigroup}) by our results, using this theorem.

\paragraph{Semigroups of Schur multipliers.}

Let $(T_t)_{t\geq 0}$ a w*-semigroup of selfadjoint contractive
Schur multipliers on $B\big(\ell^2_I\big)$. If $1 \leq p < \infty$
and $t\geq 0$, the map
$B\big(\ell^2_I\big)\xra{T_t}B\big(\ell^2_I\big)$ induces a
contraction $S^p_I \xra{T_t} S^p_I$. Using \cite[Remark 5.2]{JMX},
it is easy to see that we obtain a $C_0$-semigroup of contractions
$S^p_I \xra{T_t} S^p_I$ which satisfies inequality (\ref{inequality
NC semigroup}).

\paragraph{Semigroups of Fourier multipliers on an amenable group.}
Let $G$ be an amenable group. Let $(T_t)_{t\geq 0}$ a w*-semigroup
of selfadjoint contractive Fourier multipliers on $VN(G)$. If $1
\leq p < \infty$ and $t\geq 0$, the map $VN(G) \xra{T_t}VN(G)$
induces a contraction $L^p\big(VN(G)\big)
\xra{T_t}L^p\big(VN(G)\big)$. We obtain a $C_0$-semigroup of
contractions $L^p\big(VN(G)\big)\xra{T_t}L^p\big(VN(G)\big)$ which
satisfies inequality (\ref{inequality NC semigroup}).

\paragraph{Noncommutative Poisson semigroup.}

Let $n\geq 1$ be an integer. Recall that $\mathbb{F}_n$ denotes a
free group with $n$ generators denoted by $g_1,\ldots,g_n$. A
semigroup on $L^p\big(VN(\mathbb{F}_n)\big)$ induced by a
$w^*$-semigroup of selfadjoint completely positive unital Fourier
multipliers on $VN(\mathbb{F}_n)$ satisfies inequality
(\ref{inequality NC semigroup}). An example is provided by the
following semigroup. Any $g\in \mathbb{F}_n$ has a unique
decomposition of the form
$$
g=g_{i_1}^{k_1}g_{i_2}^{k_2}\cdots g_{i_l}^{k_l},
$$
where $l \geq 0$ is an integer, each $i_j$ belongs to
$\{1,\ldots,n\}$, each $k_j$ is a non-zero integer, and $i_j \not=
i_{j+1}$ if $1\leq j \leq l-1$. The case when $l=0$ corresponds to
the unit element $g=e_{\mathbb{F}_n}$. By definition, the length of
$g$ is defined as
$$
|g|=|k_1|+\cdots+|k_l|.
$$
This is the number of factors in the above decomposition of $g$. For
any nonnegative real number $t \geq 0$, we have a normal unital
completely positive selfadjoint map
$$
\begin{array}{cccc}
   T_t :&   \VN(\mathbb{F}_n)  &  \longrightarrow   &  \VN(\mathbb{F}_n)  \\
    &   \lambda(g)  &  \longmapsto       &  e^{-t|g|}\lambda(g).  \\
\end{array}
$$
These maps define a w*-semigroup $(T_t)_{t \geq 0} $ called the
noncommutative Poisson semigroup (see \cite{JMX} for more
information). If $1 \leq p < \infty$, this semigroup defines a
$C_0$-semigroup of contractions
$L^p\big(\VN(\mathbb{F}_n)\big)\xra{T_t}L^p\big(\VN(\mathbb{F}_n)\big)$
which satisfies inequality (\ref{inequality NC semigroup}).

\paragraph{$q$-Ornstein-Uhlenbeck semigroup.}

Suppose $-1 \leq q <1$. Let $H$ be a real Hilbert space and let
$(a_t)_{t\geq 0}$ be a $C_0$-semigroup of contractions on $H$. For
any $t\geq 0$, let $T_t=\Gamma_q(a_t)$. Then $(T_t)_{\geq 0}$ is a
w*-semigroup of normal unital completely positive maps on the von
Neumann algebra $\Gamma_q(H)$. If $1 \leq p < \infty$, this
semigroup defines a $C_0$-semigroup of contractions
$L^p\big(\Gamma_q(H)\big) \xra{T_t} L^p\big(\Gamma_q(H)\big)$ (see
\cite{JMX} for more information). This semigroup satisfies
inequality (\ref{inequality NC semigroup}).

In the case where $a_t=e^{-t}I_H$, the semigroup $(T_t)_{\geq 0}$ is
the so-called $q$-Ornstein-Uhlenbeck semigroup.

\paragraph{Modular semigroups.}
The $C_0$-semigroups of isometries satisfy inequality
(\ref{inequality NC semigroup}). Examples are provided by modular
automorphisms semigroups. Here we use noncommutative $L^p$-spaces of
a von Neumann algebra equipped with a distinguished normal faithful
state, constructed by Haagerup. We refer to \cite{PX}, and the
references therein, for more information on these spaces. Let $M$ be
a von Neumann algebra with QWEP equipped with a normal faithful
state $M\xra{\varphi}\C$. Let $\big(\sigma_t^\varphi\big)_{t\in \R}$
be the modular group of $\varphi$. If $1 \leq p < \infty$, it is
well known that $\big(\sigma_t^\varphi\big)_{t\geq 0}$ induces a
$C_0$-semigroup of isometries $L^p(M) \xra{\sigma_{t}^\varphi}
L^p(M)$ (see \cite{JuX}). This semigroup satisfies inequality
(\ref{inequality NC semigroup}).

\bigskip
In the light of Theorem \ref{th dilatation SpI}, it is natural to
ask for dilations of unital selfadjoint completely positive
semigroups of Schur multipliers. Actually, these semigroups admit a
description which allows us to construct a such dilation. %We also
%describe the semigroups of selfadjoint contractive Schur
%multipliers.
\begin{prop}
\label{description Schur mult cp} Suppose that $A$ is a matrix of
$\mathbb{M}_I$. For all $t\geq 0$, let $T_t$ be the unbounded Schur
multipliers on $B\big(\ell^2_I\big)$ associated with the matrix
\begin{eqnarray}\label{matrix semigroup}
% \nonumber to remove numbering (before each equation)
\Big[e^{-ta_{ij}}\Big]_{i,j \in I}.
\end{eqnarray}
Then the semigroup $(T_t)_{t\geq 0}$ extends to a semigroup of
selfadjoint unital completely positive Schur multipliers
        $B\big(\ell^2_I\big)\xra{T_t}B\big(\ell^2_I\big)$  if and only
        if there exists a Hilbert space $H$ and a family $(\alpha_i)_{i \in I}$ of elements of $H$ such
        that for all $t\geq 0$ the Schur multiplier $B\big(\ell^2_I\big)\xra{T_t}B\big(\ell^2_I\big)$ is associated with the
        matrix
        $$
        \Big[e^{-t\norm{\alpha_i-\alpha_j}_{H}^2}\Big]_{i,j \in I}.
        $$

%\begin{enumerate}
%  \item
%  \item The semigroup $(T_t)_{t\geq 0}$ extends to a semigroup of selfadjoint contractive Schur multipliers $B\big(\ell^2_I\big) \xra{T_t} B\big(\ell^2_I\big)$
%        if and only if there exists a Hilbert space $H$ and
%        two families $(\alpha_i)_{i \in I}$ and $(\beta_j)_{j \in I}$ of elements of $H$ such
%        that for all $t\geq 0$ the Schur multiplier $B\big(\ell^2_I\big) \xra{T_t} B\big(\ell^2_I\big)$ is associated with the
%        matrix
%        $$
%        \Big[e^{-t\norm{\alpha_i-\beta_j}_{H}^2}\Big]_{i,j \in I}.
%        $$
%\end{enumerate}
%In both case,
In this case, the Hilbert space may be chosen as a real Hilbert
space. Moreover, $(T_t)_{t\geq 0}$ is a w*-semigroup.
\end{prop}

\begin{preuve}
Now say that each $T_t$ is a selfadjoint unital completely positive
contraction means that for all $t>0$, the matrix (\ref{matrix
semigroup}) defines a real-valued positive definite kernel on
$I\times I$ in the sense of \cite[Chapter 3, Definition 1.1]{BCR}
such that for all $i\in I$ we have $a_{ii}=0$. Now, the theorem of
Schoenberg \cite[Theorem 2.2]{BCR} affirms that if $\psi$ is a
kernel then $e^{-t\psi}$ is a positive definite kernel for all $t>0$
if and only if $\psi$ is a negative definite kernel. Consequently,
the last assertion is equivalent to the fact that $A$ defines a
real-valued negative definite kernel which vanishes on the diagonal
of $I \times I$. Finally, the characterization of real-valued
definite negative kernel of \cite[Proposition 3.2]{BCR} gives the
equivalence with the required description.

The assertion concerning the choice of the Hilbert space is clear.
Finally, using \cite[Remark 5.2]{JMX}, it is easy to see that
$(T_t)_{t\geq 0}$ is a w*-semigroup.
\end{preuve}

The next proposition is inspired by the work \cite{JuX}.

\begin{prop}
Let $(T_t)_{t\geq 0}$ be a w*-semigroup of selfadjoint unital
completely positive Schur multipliers on $B\big(\ell^2_I\big)$.
Then, there exists a hyperfinite von Neumann algebra $M$ equipped
with a semifinite normal faithful trace, a w*-semigroup $(U_t)_{t
\geq 0}$ of unital trace preserving *-auto\-mor\-phisms of $M$, a
unital trace preserving one-to-one normal $*$-homomorphism
$B\big(\ell^2_I\big) \xra{J} M$ such that
\begin{eqnarray*}
T_t=\mathbb{E}U_tJ.
\end{eqnarray*}
for any $t \geq 0$, where $M \xra{\mathbb{E}} B\big(\ell^2_I\big)$
is the canonical faithful normal trace preserving conditional
expectation associated with $J$.
\end{prop}

\begin{preuve}
By Proposition \ref{description Schur mult cp}, there exists a real
Hilbert space $H$ and a family $(\alpha_j)_{j \in I}$ of elements of
$H$ such that, for all $t\geq 0$, the Schur multiplier
$B\big(\ell^2_I\big)\xra{T_t}B\big(\ell^2_I\big)$ is associated with
the matrix
$$
\Big[e^{-t\norm{\alpha_j-\alpha_k}_{H}^2}\Big]_{j,k \in I}.
$$
Let $\mu$ be a gaussian measure on $H$, i.e. a probability space
$(\Omega,\mu)$ together with a measurable function $\Omega\xra{w}H$
such that, for all $h\in H$, we have
$$
e^{-\norm{h}_{H}^2} = \int_{\Omega} e^{i\langle
h,w(\omega)\rangle_{H}}d\mu(\omega)
$$
where $i^2=-1$. We define the von Neumann algebra
$M=L^\infty(\Omega)\otvn B\big(\ell^2_I\big)$. Note that $M$ is a
hyperfinite von Neumann algebra. We equip the von Neumann algebra
$M$ with the faithful semifinite normal trace
$\tau_M=\int_{\Omega}\cdot\ d\mu\ot \tr$. Note that, by
\cite[Theorem 1.22.13]{Sak}, we have a $*$-isomorphism
$M=L^\infty\big(\Omega,B(\ell^2_I)\big)$. We define the canonical
injective normal unital $*$-homomorphism
$$
\begin{array}{cccc}
  J:  &  B\big(\ell^2_I\big)   &  \longrightarrow   &   L^\infty(\Omega)\otvn B\big(\ell^2_I\big) \\
      &       x        &  \longmapsto       &  1 \ot x.  \\
\end{array}
$$
It is clear that the map $J$ preserves the traces. We denote by $M
\xra{\mathbb{E}} B\big(\ell^2_I\big)$ the canonical faithful normal
trace preserving conditional expectation of $M$ onto
$B\big(\ell^2_I\big)$. For all $\omega\in \Omega$ and $t>0$ let
$D_t(\omega)$ be the diagonal matrix of $B\big(\ell^2_I\big)$
defined by
$$
D_t(\omega)=\bigg[\delta_{j,k}e^{i \sqrt{t}\langle
\alpha_j,w(\omega)\rangle_{\ell^2_I}}\bigg]_{j,k \in I}.
$$
Note that, for all $t>0$, the map $\Omega\xra{D_t}B(H)$ defines an
unitary element of $L^\infty\big(\Omega,B(\ell^2_I)\big)$. Now, for
all $t\geq 0$ we define the linear map
$$
\begin{array}{cccc}
 U_t:   &  L^\infty\big(\Omega,B(\ell^2_I)\big)   &  \longrightarrow   &  L^\infty\big(\Omega,B(\ell^2_I)\big)  \\
        &  f   &  \longmapsto       &  D_t fD_t^*. \\
\end{array}
$$
If $t\geq 0$, it is easy to see that the map $U_t$ is a trace
preserving $*$-automorphism of $M$. For all $x\in
B\big(\ell^2_I\big)$, we have
\begin{eqnarray*}
% \nonumber to remove numbering (before each equation)
 \mathbb{E}U_tJ(x)  &=&  \mathbb{E}U_t(1\ot x)  \\
   &=&  \int_{\Omega}  D_t(\omega) (1\ot x)D_t(\omega)^*d\mu(\omega) \\
   &=&  \int_{\Omega}\bigg[e^{i \sqrt{t}\langle \alpha_j-\alpha_k,w(\omega)\rangle_{H}}x_{jk}\bigg]_{j,k\in I}d\mu(\omega)\\
   &=&  \Big[e^{-t\norm{\alpha_j-\alpha_k}_{H}^2}x_{jk}\Big]_{j,k \in I}\\
   &=&  T_t(x).
\end{eqnarray*}
Thus, for all $t\geq 0$, we have
$$
T_t=\mathbb{E}U_tJ.
$$
The assertion concerning the regularity of the semigroup is easy and
left to the reader.
\end{preuve}

\vspace{0.5cm}

In the same vein, it is not difficult to construct a dilation of the
noncommutative Poisson semigroup. The result was already known to F.
Lust-Piquard. Moreover, it is easy to dilate the $C_0$-semigroups of
contractions $L^p\big(\Gamma_q(H)\big) \xra{\Gamma_q(a_t)}
L^p\big(\Gamma_q(H)\big)$, with \cite[Theorem 8.1]{SNF}.
\vspace{0.7cm}

Finally, we have the next result analogue to Corollary
\ref{dilationQWEP}. One can prove this proposition with a similar
argument.
\begin{prop}
\label{dilation semigroup QWEP} Suppose $1< p < \infty$. Let
$(T_t)_{t\geq 0}$ be a $C_0$-semigroup of contractions on a
noncommutative $L^p$-space $L^p(M)$. Suppose that there exists a
noncommutative $L^p$-space $L^p(N)$ where $N$ has QWEP, a
$C_0$-semigroup $(U_t)_{t \geq 0}$ of isometric operators on
$L^p(N)$, an isometric embedding $L^p(M)\xra{J}L^p(N)$ and a
contractive map $L^p(N)\xra{Q}L^p(M)$ such that,
$$
T_t=QU_tJ.
$$
for any $t\geq 0$. Then, for all $b \in L^1(\mathbb{R})$ with
support in ${\mathbb{R}^+}$, we have the estimate
$$
\Bgnorm{\int_{0}^{+\infty}b(t)T_t dt}_{L^p(M) \xra{} L^p(M)} \leq
\Bgnorm{\int_{0}^{+\infty}b(t)\mathcal{T}_t dt}_{cb,L^p(\mathbb{R})
\xra{}L^p(\mathbb{R})}.
$$
Moreover, if $L^p(N)$ is a commutative $L^p$-space $L^p(\Omega)$, we
have, for all $b \in L^1(\mathbb{R})$ with support in
${\mathbb{R}^+}$, the estimate
$$
\Bgnorm{\int_{0}^{+\infty}b(t)T_t dt}_{L^p(\Omega) \xra{}
L^p(\Omega)} \leq \Bgnorm{\int_{0}^{+\infty}b(t)\mathcal{T}_t
dt}_{L^p(\mathbb{R}) \xra{}L^p(\mathbb{R})}.
$$
\end{prop}

\vspace{0.5cm}

This proposition allows us to give alternate proofs for some results
of this section. By example, using \cite[Theorem 8.1]{SNF} of
dilation of $C_0$-semigroups on a Hilbert space, we deduce that the
both inequalities (\ref{inequality semigroup}) and (\ref{inequality
NC semigroup}) are true for $p=2$. By using \cite{Fen}, we see that
the $C_0$-semigroups of operators which admit a contractive majorant
satisfy inequality (\ref{inequality semigroup}), for $1 < p <
\infty$ .

%%%%%%%%%%%%%%%%%%%%%%%%%%%%%%%%%%%%%%%%%%%%%%%%%%%%%%%%%%%%%%%%
%%%%%%%%%%%%%%%%%%%%%%%%%%%%%%%%%%%%%%%%%%%%%%%%%%%%%%%%%%%%%%%%

\section{Final remarks}

%%%%%%%%%%%%%%%%%%%%%%%%%%%%%%%%%%%%%%%%%%%%%%%%%%%%%%%%%%%%%%%%
%%%%%%%%%%%%%%%%%%%%%%%%%%%%%%%%%%%%%%%%%%%%%%%%%%%%%%%%%%%%%%%%

We begin by observing that the inequalities (\ref{Matsaev}) and
(\ref{LpMatsaev}) are true for any complex polynomial $P$ of degree
1 and any contraction $T$. Indeed, suppose that $P(z)=a+bz$, then it
is easy to see that $\norm{P}_{2}=|a|+|b|$. Thus, for all $1 \leq p
\leq \infty$, we have $\norm{P}_{p}=\norm{P}_{p,S^p}=|a|+|b|$.

Now we will determine the real polynomials of higher degree with a
similar property.
\begin{prop}
Let $\displaystyle P=\sum_{k=0}^n a_{k}z^k$ be a real polynomial
such that $a_k\not=0$ for any $0\leq k\leq n$. The following
assertions are equivalent.
\begin{enumerate}
  \item For all $1 < p <\infty$,  we have $\displaystyle \norm{P}_{p}=\sum_{k=0}^n
  |a_{k}|$.
  \item For all $1 < p <\infty$, we have $\displaystyle  \norm{P}_{p,S^p}=\sum_{k=0}^n
  |a_{k}|$.
  \item There exists $1 < p <\infty$ such that $\displaystyle \norm{P}_{p}=\sum_{k=0}^n
  |a_{k}|$.
  \item There exists $1 < p <\infty$ such that $\displaystyle  \norm{P}_{p,S^p}=\sum_{k=0}^n
  |a_{k}|$.
  \item The coefficients $a_k$ have the same sign or the signs
  of the $a_k$ are alternating (i.e. for any integer $0 \leq k \leq n-1$ we have $a_ka_{k+1}\leq 0$).
\end{enumerate}
In this case, for the polynomial $P$ and any contraction $T$, the
inequalities (\ref{Matsaev}) and (\ref{LpMatsaev}) are true.
\end{prop}

\begin{preuve}
First we will show that  $\norm{P}_{2}=\sum_{k=0}^n |a_{k}|$ is
equivalent to the last assertion. Recall that
$\norm{P}_{2}=\sup_{|z|=1}|P(z)|$. On the one hand, for all $0 \leq
\theta \leq 2\pi$, we have
\begin{eqnarray*}
% \nonumber to remove numbering (before each equation)
 \lefteqn{ \left|\sum_{k=0}^n a_ke^{ki\theta}\right|^2
 =\Bigg(\sum_{k=0}^n a_k \cos(k\theta)\Bigg)^2+\Bigg(\sum_{k=0}^n a_k
 \sin(k\theta)\Bigg)^2}\\
   &=& \sum_{k=0}^n a_{k}^2 \cos^2(k\theta)+2\sum_{0\leq k<l\leq n}a_{k}a_{l}\cos(k\theta)\cos(l\theta)+\sum_{k=0}^n a_{k}^2 \sin^2(k\theta)\\
   &&+2\sum_{0\leq k<l\leq n}a_{k}a_{l}\sin(k\theta)\sin(l\theta) \\
   &=& \sum_{k=0}^n a_{k}^2 +2\sum_{0\leq k<l\leq n}a_{k}a_{l}\cos\big((k-l)\theta\big).
\end{eqnarray*}
On the other hand, we have the equality
$$
\Bigg(\sum_{k=0}^n |a_{k}|\Bigg)^2= \sum_{k=0}^n
a_{k}^2+2\sum_{0\leq k<l\leq n}|a_{k}a_{l}|.
$$
Then  $P$ satisfies $\norm{P}_{2}=\sum_{k=0}^n |a_{k}|$ if and only
if
$$
\sum_{0\leq k<l\leq
n}a_{k}a_{l}\cos\big((k-l)\theta\big)=\sum_{0\leq k<l\leq
n}|a_{k}a_{l}|.
$$
This last assertion means that for all $0\leq k<l\leq n$ we have $
\cos\big((k-l)\theta\big)={\rm sign} (a_ka_l)$. It is easy to see
that this last assertion is equivalent to the assertion 5.

Now, it is trivial that the equality $\norm{P}_{2}=\sum_{k=0}^n
|a_{k}|$ implies the assertions 1 and 2, that 1 implies 3 and that 2
implies 4. Now we show that 4 implies $\norm{P}_{2}=\sum_{k=0}^n
|a_{k}|$. By interpolation, we have
$$
\norm{P}_{\infty,S_\infty}=\sum_{k=0}^n |a_{k}|
=\norm{P}_{p,S^p}\leq
\big(\norm{P}_{\infty,S_\infty}\big)^{1-\frac{2}{p}}\big(
\norm{P}_{2,S_2}\big)^{\frac{2}{p}}.
$$
Moreover, it is easy to see that $\norm{P}_{2}=\norm{P}_{2,S_2}$.
Then we obtain
$$
\big(\norm{P}_{\infty,S_\infty}\big)^{\frac{2}{p}} \leq
\big(\norm{P}_{2,S_2}\big)^{\frac{2}{p}} =
\big(\norm{P}_{2}\big)^{\frac{2}{p}}.
$$
And finally we have
$$
\sum_{k=0}^n |a_{k}| = \norm{P}_{\infty,S_\infty} \leq \norm{P}_{2}.
$$
The proof that the assertion 3 implies $\norm{P}_{2}=\sum_{k=0}^n
|a_{k}|$ is similar.
\end{preuve}

\bigskip

%%%%%%%%%%%%%%%%%%%%%%%%%%%%%%%%%%%%%%%%%%%%%%%%%%%%%%%%%%%%%%%%%%%%%%%%%%%%%%%%%%%%%%%%%%%%%%%%%%%%%%%%%%%%%%%%%%%%%%%%%%%%%%

\textbf{Acknowledgement}. I wish to thank my thesis advisor
Christian Le Merdy for his support and advice, Eric Ricard for
fruitful discussions and Jesse Peterson for suggest me the proof of
Proposition \ref{prop AVN has QWEP}.

%%%%%%%%%%%%%%%%%%%%%%%%%%%%%%%%%%%%%%%%%%%%%%%%%%%%%%%%%%%%%%%%%%%%%%%%%%%%%%%%%%%%%%%%%%%%%%%%%%%%%%%%%%%%%%%%%%%%%%%%%%%%%

\bigskip\footnotesize{
\n Laboratoire de Math\'ematiques, Universit\'e de Franche-Comt\'e,
25030 Besan\c{c}on Cedex,  France\\
cedric.arhancet@univ-fcomte.fr\hskip.3cm


\small
\begin{thebibliography}{79}

\bibitem[AbA]{AbA}
Y. Abramovich and C. Aliprantis.
\newblock {\em An invitation to operator theory}.
\newblock American Mathematical Society, Providence, 2002.

\bibitem[AkS]{AkS}
M. Akcoglu and L. Sucheston.
\newblock Dilations of positive contractions on $L_p$ spaces.
\newblock  {\em Canad. Math. Bull.} 20, no. 3, 285--292, 1977.

\bibitem[ALM]{ALM}
C. Arhancet and C. Le Merdy.
\newblock Dilation of Ritt operators on $L^p$-spaces.
\newblock Preprint, arXiv:1106.1513.

\bibitem[Arh]{Arh}
C. Arhancet.
\newblock Unconditionality, Fourier multipliers and Schur multipliers.
\newblock  {\em Colloq. Math.} 127, 17--37, 2012.

\bibitem[BCR]{BCR}
C. Berg, J. Christensen and P. Ressel.
\newblock {\em Harmonic analysis on semigroups. Theory of positive definite and related functions}.
\newblock Springer-Verlag, New York, 1984.

\bibitem[BDJ]{BDJ}
N. Brown, K. Dykema and K. Jung.
\newblock Free entropy dimension in amalgamated free products. With an appendix by Wolfgang Lück.
\newblock {\em  Proc. Lond. Math. Soc.} 97, no. 2, 339--367, 2008.

\bibitem[BeL]{BeL}
J. Bergh and J. L{\"o}fstr{\"o}m.
\newblock {\em Interpolation spaces}.
\newblock Springer-Verlag, Berlin, 1976.

\bibitem[BGM]{BGM}
E. Berkson, T. Gillespie and P. Muhly.
\newblock Generalized analyticity in UMD spaces.
\newblock {\em Ark. Mat.} 27, no. 1: 1--14, 1989.

\bibitem[BiS]{BiS}
M. Birman and M. Solomyak.
\newblock Double operator integrals in a Hilbert space.
\newblock {\em Integral Equations Operator Theory} 47, no. 2: 131--168, 2003.

\bibitem[BKS]{BKS}
M. Bo\.{z}ejko, B. Kümmerer and R. Speicher.
\newblock $q$-Gaussian processes: non-commutative and classical aspects.
\newblock {\em  Comm. Math. Phys.} 185, no. 1, 129--154, 1997.

\bibitem[Bla]{Bla}
B. Blackadar.
\newblock {\em Operator algebras. Theory of C*-algebras and von Neumann algebras}.
\newblock Springer-Verlag, Berlin, 2006.

\bibitem[BlD]{BlD}
E. Blanchard and K. Dykema.
\newblock Embeddings of reduced free products of operator algebras..
\newblock {\em Pacific J. Math.} 199, no. 1: 1--19, 2001.

\bibitem[BoS]{BoS}
M. Bo\.{z}ejko and R. Speicher.
\newblock Completely positive maps on Coxeter groups, deformed commutation relations, and operator spaces.
\newblock {\em Math. Ann.} 300, no. 1: 97--120, 1994.

\bibitem[Boz]{Boz}
M. Bo\.{z}ejko.
\newblock Bessis-Moussa-Villani conjecture and generalized Gaussian random variables.
\newblock {\em Infin. Dimens. Anal. Quantum Probab. Relat. Top.} 11, no. 3: 313--321, 2008.

\bibitem[Cha]{Cha}
Chatterji, S. D.
\newblock  Martingale convergence and the Radon-Nikodym
theorem in Banach spaces.
\newblock  {\em Math. Scand.} 22, 21--41, 1968.

\bibitem[Con]{Con}
A. Connes.
\newblock Classification of injective factors.
\newblock {\em Ann. of Math.} 104, no. 1: 73--115, 1976.

\bibitem[CoW]{CoW}
R. Coifman and G. Weiss.
\newblock {\em Transference methods in analysis.}
\newblock Conference Board of the Mathematical Sciences Regional Conference Series in Mathematics, No.
31. American Mathematical Society, Providence, R.I., 1976.

\bibitem[DCH]{DCH}
J. De Cannière and U. Haagerup.
\newblock Multipliers of the Fourier algebras of some simple Lie groups and their discrete subgroups.
\newblock {\em Amer. J. Math.} 107, no. 2: 455--500, 1985.

\bibitem[DeL]{DeL}
K. De Leeuw.
\newblock On $L_p$ multipliers.
\newblock {\em Ann. of Math.} 81: 364--379, 1965.

\bibitem[DiU]{DU}
J. Diestel and J. Uhl.
\newblock {\em Vector measures}.
\newblock American Mathematical Society, Providence, 1977.

\bibitem[Dr]{Dr}
S. W. Drury.
\newblock A counterexample to a conjecture of Matsaev.
\newblock {\it Lin. Alg. and its Appl.} 435, 323-329 (2011).

\bibitem[EfR]{ER1}
E. Effros and Z-J. Ruan.
\newblock {\em Operator spaces}.
\newblock Oxford University Press, 2000.

\bibitem[EfP]{EfP}
E. Effros and M. Popa.
\newblock Feynman diagrams and Wick products associated with q-Fock space.
\newblock {\em Proc. Natl. Acad. Sci. USA} 100, no. 15: 8629--8633, 2003.

\bibitem[Fen]{Fen}
G. Fendler.
\newblock Dilations of one parameter semigroups of positive contractions on $L^p$ spaces.
\newblock {\em Canad. J. Math.} 49, no. 4: 736--748, 1997.

\bibitem[Jod]{Jod}
M. Jodeit.
\newblock Restrictions and extensions of Fourier multipliers.
\newblock {\em Studia Math.} 34: 215--226, 1970.

\bibitem[JMX]{JMX}
M. Junge,  C. Le Merdy and Q. Xu.
\newblock {\em $H^\infty$ functional calculus and square functions on
noncommutative $L^p$-spaces.}
\newblock  Astérisque,  no. 305, 2006.

\bibitem[JLM]{JLM}
M. Junge and C. Le Merdy.
\newblock Dilations and rigid factorisations on noncommutative $L_p$-spaces.
\newblock {\em J. Funct. Anal.} 249: 220--252, 2007.

\bibitem[JuX]{JuX}
M. Junge and Q. Xu.
\newblock Noncommutative maximal ergodic theorems.
\newblock {\em J. Amer. Math. Soc.} 20, no. 2: 385--439, 2007.

\bibitem[Jun]{Jun}
M. Junge.
\newblock Fubini's theorem for ultraproducts of noncommmutative $L_p$-spaces II.
\newblock Preprint.

\bibitem[Kir]{Kir}
E. Kirchberg.
\newblock On nonsemisplit extensions, tensor products and exactness of group C*-algebras.
\newblock {\em Invent. Math.} 112, no. 3: 449--489, 1993.

\bibitem[Kit]{Kit}
A. Kitover.
\newblock A question in connection with Matsaev's conjecture.
\newblock {\em Linear and Complex Analysis Problem Book 3, Part I}, Lecture Notes in Math. 1573: 247, Springer, Berlin, 1985.

\bibitem[Lar]{Lar}
R. Larsen.
\newblock {\em An introduction to the theory of multipliers.}
\newblock Springer-Verlag, 1971.

\bibitem[LaS]{LaS}
V. Lafforgue and M. de la Salle.
\newblock Noncommutative $L^p$-spaces without the completely bounded approximation property.
\newblock {\em Duke Math. J.} 160 , no. 1, 71--116, 2011.

\bibitem[MeN]{MeN}
P. Meyer-Nieberg.
\newblock {\em Banach lattices}.
\newblock Springer-Verlag, Berlin, 1991.

\bibitem[NeR]{NeR}
S. Neuwirth and E. Ricard.
\newblock Transfer of
Fourier multipliers into Schur multipliers and sumsets in a discrete
group.
\newblock {\em Canad. J. Math.} 63, no. 5, 1161--1187, 2011.  %1001.5332

\bibitem[Nik1]{Nik1}
N. K. Nikolski.
\newblock Five problems on invariant subspaces.
\newblock {\em J. Sov. Math.}, 2: 441--450, 1974.

\bibitem[Nik2]{Nik2}
N. Nikolski.
\newblock {\em Operators, functions, and systems: an easy reading, vol 2.}
\newblock American Mathematical Society, 2002.

\bibitem[Nou]{Nou}
A. Nou.
\newblock Asymptotic matricial models and QWEP property for $q$--Araki--Woods algebras.
\newblock {\em J. Funct. Anal.} 232, no. 2: 295--327, 2006.

\bibitem[Oza]{Oza}
N. Ozawa.
\newblock About the QWEP conjecture.
\newblock {\em Internat. J. Math.} 15, no. 5: 501--530, 2004.

\bibitem[Pau]{Pau}
V. Paulsen.
\newblock {\em Completely bounded maps and operator algebras}.
\newblock Cambridge Univ. Press, 2002.

\bibitem[Pel1]{Pel1}
V.V. Peller.
\newblock An Analogue of an inequality of J. von Neumann, isometric dilation of contractions, and
approximation by isometries in spaces of measurable functions.
\newblock {\em Proc. Steklov Inst. Math.} 1 : 101--145, 1983.

\bibitem[Pel2]{Pel2}
V.V. Peller.
\newblock Estimates of operator polynomials on the Schaten - von Neumann classes.
\newblock {\em Linear and Complex Analysis Problem Book 3, Part I}, Lecture Notes in Math. 1573: 244--246, Springer, Berlin, 1985.

\bibitem[Pis1]{Pis222}
G. Pisier.
\newblock Regular operators between non-commutative $L_p$-spaces.
\newblock {\em Bull. Sci. Math.} 119 : 95--118, 1995.

\bibitem[Pis2]{Pis3}
G. Pisier.
\newblock Non-commutative vector valued {$L\sb p$}-spaces and completely {$p$}-summing maps.
\newblock {\em Astérisque}, 247, 1998.

\bibitem[Pis3]{Pis33}
G. Pisier.
\newblock {\em Similarity problems and completely bounded maps}, volume 1618 of
  {\em Lecture Notes in Mathematics}.
\newblock Springer-Verlag, expanded edition, 2001.

\bibitem[Pis4]{Pis4}
G. Pisier.
\newblock {\em Introduction to operator space theory}.
\newblock Cambridge University Press, Cambridge, 2003.

\bibitem[PiX]{PX}
G. Pisier and Q. Xu.
\newblock Non-commutative $L_p$-spaces.
\newblock {\em Handbook of the Geometry of Banach Spaces volume II}: 1459--1517, 2003.

\bibitem[Ric]{Ric}
E. Ricard.
\newblock A Markov dilation for self-adjoint Schur multipliers.
\newblock {\em Proc. Amer. Math. Soc.} 136: 4365--4372, 2008.

\bibitem[Sak]{Sak}
S. Sakai.
\newblock {\em C*-algebras and W*-algebras}.
\newblock Springer-Verlag, New York, 1971.

\bibitem[Sun]{Sun}
V. Sunder.
\newblock {\em An invitation to von Neumann algebras}.
\newblock Springer-Verlag, New York, 1987.

\bibitem[SNF]{SNF}
B. Sz.-Nagy and C. Foias.
\newblock {\em Harmonic analysis of operators on Hilbert space}.
\newblock North-Holland, Amsterdam-London, 1970.

\bibitem[Str]{Stra}
S. Stratila.
\newblock {\em Modular theory in operator algebras}.
\newblock Taylor and Francis, 1981.

\bibitem[Tak3]{Tak3}
M. Takesaki.
\newblock {\em Theory of operator algebras. III.}
\newblock Springer-Verlag, Berlin, 2003.

\bibitem[Ued]{Ued}
Y. Ueda.
\newblock Amalgamated free product over Cartan subalgebra.
\newblock {\em Pacific J. Math.} 191, no. 2: 359--392, 1999.

\end{thebibliography}
\end{document}